
%
\magnification1200
\pretolerance=100
\tolerance=200
\hbadness=1000
\vbadness=1000
\linepenalty=10
\hyphenpenalty=50
\exhyphenpenalty=50
\binoppenalty=700
\relpenalty=500
\clubpenalty=5000
\widowpenalty=5000
\displaywidowpenalty=50
\brokenpenalty=100
\predisplaypenalty=7000
\postdisplaypenalty=0
\interlinepenalty=10
\doublehyphendemerits=10000
\finalhyphendemerits=10000
\adjdemerits=160000
\uchyph=1
\delimiterfactor=901
\hfuzz=0.1pt
\vfuzz=0.1pt
\overfullrule=5pt
\hsize=146 true mm
\vsize=8.9 true in
\maxdepth=4pt
\delimitershortfall=.5pt
\nulldelimiterspace=1.2pt
\scriptspace=.5pt
\normallineskiplimit=.5pt
\mathsurround=0pt
\parindent=20pt
\catcode`\_=11
\catcode`\_=8
\normalbaselineskip=12pt
\normallineskip=1pt plus .5 pt minus .5 pt
\parskip=6pt plus 3pt minus 3pt
\abovedisplayskip = 12pt plus 5pt minus 5pt
\abovedisplayshortskip = 1pt plus 4pt
\belowdisplayskip = 12pt plus 5pt minus 5pt
\belowdisplayshortskip = 7pt plus 5pt
\normalbaselines
\smallskipamount=\parskip
 \medskipamount=2\parskip
 \bigskipamount=3\parskip
\jot=3pt
%
%
\def\ref#1{\par\noindent\hangindent2\parindent
 \hbox to 2\parindent{#1\hfil}\ignorespaces}
%
%
\font\tenss=cmss10
\font\sevenss=cmss8 at 7pt
\font\fivess=cmss8 at 5pt
\newfam\ssfam %
\textfont\ssfam=\tenss
\scriptfont\ssfam=\sevenss
\scriptscriptfont\ssfam=\fivess
%
%
%
%
%
%
%
%
%
\catcode`\_=11
\def\suf_fix{}
\def\scaled_rm_box#1{%
 \relax
 \ifmmode
   \mathchoice
    {\hbox{\tenrm #1}}%
    {\hbox{\tenrm #1}}%
    {\hbox{\sevenrm #1}}%
    {\hbox{\fiverm #1}}%
 \else
  \hbox{\tenrm #1}%
 \fi}
\def\suf_fix_def#1#2{\expandafter\def\csname#1\suf_fix\endcsname{#2}}
\def\I_Buchstabe#1#2#3{%
 \suf_fix_def{#1}{\scaled_rm_box{I\hskip-0.#2#3em #1}}
}
\def\rule_Buchstabe#1#2#3#4{%
 \suf_fix_def{#1}{%
  \scaled_rm_box{%
   \hbox{%
    #1%
    \hskip-0.#2em%
    \lower-0.#3ex\hbox{\vrule height1.#4ex width0.07em }%
   }%
   \hskip0.50em%
  }%
 }%
}
\I_Buchstabe B22
\rule_Buchstabe C51{34}
\I_Buchstabe D22
\I_Buchstabe E22
\I_Buchstabe F22
\rule_Buchstabe G{525}{081}4
\I_Buchstabe H22
\I_Buchstabe I20
\I_Buchstabe K22
\I_Buchstabe L20
\I_Buchstabe M{20em }{I\hskip-0.35}
\I_Buchstabe N{20em }{I\hskip-0.35}
\rule_Buchstabe O{525}{095}{45}
\I_Buchstabe P20
\rule_Buchstabe Q{525}{097}{47}
\I_Buchstabe R21 
\rule_Buchstabe U{45}{02}{54}
\suf_fix_def{Z}{\scaled_rm_box{Z\hskip-0.38em Z}}
\catcode`\"=12
\newcount\math_char_code
\def\suf_fix_math_chars_def#1{%
 \ifcat#1A
  \expandafter\math_char_code\expandafter=\suf_fix_fam
  \multiply\math_char_code by 256
  \advance\math_char_code by `#1
  \expandafter\mathchardef\csname#1\suf_fix\endcsname=\math_char_code
  \let\next=\suf_fix_math_chars_def
 \else
  \let\next=\relax
 \fi
 \next}
%
%
%
%
\def\font_fam_suf_fix#1#2 #3 {%
 \def\suf_fix{#2}
 \def\suf_fix_fam{#1}
 \suf_fix_math_chars_def #3.
}
\font_fam_suf_fix
 0rm
 ABCDEFGHIJKLMNOPQRSTUVWXYZabcdefghijklmnopqrstuvwxyz
\font_fam_suf_fix
 2scr
 ABCDEFGHIJKLMNOPQRSTUVWXYZ
\font_fam_suf_fix
 \slfam sl
 ABCDEFGHIJKLMNOPQRSTUVWXYZabcdefghijklmnopqrstuvwxyz
\font_fam_suf_fix
 \bffam bf
 ABCDEFGHIJKLMNOPQRSTUVWXYZabcdefghijklmnopqrstuvwxyz
\font_fam_suf_fix
 \ttfam tt
 ABCDEFGHIJKLMNOPQRSTUVWXYZabcdefghijklmnopqrstuvwxyz
\font_fam_suf_fix
 \ssfam
 ss
 ABCDEFGHIJKLMNOPQRSTUVWXYZabcdefgijklmnopqrstuwxyz
\catcode`\_=8
\def\Cdss{{\fam\ssfam
    \mkern 4.2 mu \mathchoice%
    {\vrule height 6.5pt depth -.55pt width 1pt}%
    {\vrule height 6.5pt depth -.57pt width 1pt}%
    {\vrule height 4.55pt depth -.28pt width .8pt}%
    {\vrule height 3.25pt depth -.19pt width .6pt}%
    \mkern -6.3mu C}}%
\def\Ndss{{\fam\ssfam I\mkern -2.5mu N}}%
\def\Qdss{{\fam\ssfam
    \mkern 3.8 mu \mathchoice%
    {\vrule height 6.5pt depth -.67pt width 1pt}%
    {\vrule height 6.5pt depth -.7pt width 1pt}%
    {\vrule height 4.55pt depth -.44pt width .7pt}%
    {\vrule height 3.25pt depth -.3pt width .5pt}%
    \mkern -5.9mu Q}}%
\def\Zdss{{\fam\ssfam Z\mkern-8.1mu Z}}%
%
%
%
%
\font\teneuf=eufm10 
\font\seveneuf=eufm7
\font\fiveeuf=eufm5
\newfam\euffam \def\euf{\fam\euffam\teneuf} 
\textfont\euffam=\teneuf \scriptfont\euffam=\seveneuf
\scriptscriptfont\euffam=\fiveeuf
       \def\afr{{\euf a}}
       
       \def\cfr{{\euf c}}
       
       \def\efr{{\euf e}}
       
       \def\gfr{{\euf g}}
       \def\hfr{{\euf h}}

       \def\lfr{{\euf l}}
       \def\mfr{{\euf m}}

       \def\sfr{{\euf s}}
       
       \def\ufr{{\euf u}}

       \def\xfr{{\euf x}}
       
       \def\zfr{{\euf z}}
\parindent=0pt
\let\lra=\longrightarrow
\let\lla=\longleftarrow
\let\ul=\underline
\def\kat#1{\underline{\cal #1}}
\def\N{{\Ndss}}
\def\dZ{{\Zdss_p}}
\def\dQ{{\Qdss_p}}
\def\dC{{\Cdss_p}}
\def\cF{{\cal F}}
\def\cI{{\cal I}}
\def\cO{{\cal O}}
\def\M{{\cal M}}
\def\Ind{{\rm Ind}}
\def\Hom{\mathop{\rm Hom}\nolimits}
\def\Lie{\mathop{\rm Lie}\nolimits}
\def\L{{\cal L}}
\def\LC{\kat{LC}}
\def\Ad{\hbox{\rm Ad}}
\def\ad{\hbox{\rm ad}}
\def\to{\rightarrow}
\def\morphism#1#2#3#4{{\matrix{\hfill #1 &\lra& #2\hfill\cr\cr
\hfill #3 &\longmapsto& #4\hfill\cr} }}
\hyphenation{Ba-nach Stein-haus}
\centerline{\bf Locally analytic distributions and $p$-adic representation theory,}
\centerline{\bf with applications to $GL_2$}

\medskip

\centerline{\bf P. Schneider, J. Teitelbaum}

\bigskip

Let $L/\dQ$ be a finite extension, and let $G$ be a locally
$L$-analytic group such as the $L$-points of an algebraic group
over $L$. The categories of smooth and of finite-dimensional
algebraic representations of $G$ may be viewed as subcategories of
the category of continuous representations  of $G$ in locally
convex $L$-vector spaces. This larger category contains many
interesting new objects, such as the action of $G$ on global
sections of equivariant vector bundles over $p$-adic symmetric
spaces and other representations arising from the theory of
$p$-adic uniformization as studied for example in [ST]. A workable
theory of continuous representations of $G$ in $L$-vector spaces
offers the opportunity to unify these disparate examples in a
single theoretical framework.

There are a number of technical obstacles to developing a
reasonable theory of such representations. For example, there are
no unitary representations over $L$, and continuous, or even
locally analytic, functions on $G$ are not integrable against Haar
measure. As a result, even for compact groups one is forced to
consider representations of $G$ in fairly general locally convex
vector spaces.  In such situations one encounters a range of
pathologies apparent even in the theory of representations of real
Lie groups in Banach spaces.   For this reason one must formulate
some type of ``finiteness'' or ``admissibility'' condition in
order to have a manageable theory.

In this paper we introduce a restricted category of continuous
representations of  locally $L$-analytic groups $G$ in locally
convex $K$-vector spaces, where $K$ is a spherically complete
non-archimedean extension field of $L$. We call the objects of
this category ``admissible'' representations and we establish some
of their basic properties. Most importantly we show that (at least
when $G$ is compact) the category of admissible representations in
our sense can be algebraized; it is faithfully full
(anti)-embedded into the category of modules over the locally
analytic distribution algebra $D(G,K)$ of $G$ over $K$. We may
then replace the topological notion of irreducibility with the
algebraic property of simplicity as $D(G,K)$-module. Our hope is
that our definition of admissible representation may be used as a
foundation for a general theory of continuous $K$-valued
representations of locally $L$-analytic groups.

As  an application of our theory, we prove the topological
irreducibility of generic members of the $p$-adic principal series
of $GL_{2}(\dQ)$.   This result was claimed by Morita for general
$L$, not only $\dQ$, in [Mor] Thm. 1(i), but his method to deduce
irreducibility from something weaker that he calls local
irreducibility is flawed. We believe that this failure is in fact
an unavoidable consequence of the above mentioned pathologies of
Banach space representations. Our completely different method is
based instead on our algebraization theory and leads to the much
stronger result of the algebraic simplicity of the corresponding
$D(GL_2 (\dQ),K)$-modules. Morita also described the intertwining
operators between the various $p$-adic principal series
representations ([Mor] Thm. 2); again we use our methods to
establish a slightly stronger version of this.

We rely heavily on two comparatively old results. The first one is
Amice's Fourier isomorphism in [Am2] between $D(\dZ,K)$ and the
ring of power series over $K$ convergent on the open $p$-adic unit
disk. The second one is Lazard's result in [Laz] that in the
latter ring every divisor is principal. The structure of the ring
$D(o_L)$, for the ring of integers $o_L$ in a general finite
extension $L$ of $\dQ$, is unclear at present, and for this reason
we obtain our irreducibility theorem only for $\dQ$.

We make extensive use of results from $p$-adic functional
analysis. Where possible, we give references to both the classical
and  the $p$-adic literature.  In some cases, however, we have not
found suitable $p$-adic references and we have given references to
the needed result in the classical case. These classical results
are general properties of Fr\'{e}chet spaces and their $p$-adic
analogues are not difficult to obtain (see [NFA] for a systematic
presentation of $p$-adic methods).

Most of this work was done during a stay of the first author at
the University of Illinois at Chicago. He wants to express his
gratitude for the support during this very pleasant and fruitful
time.

\bigskip

{\bf 1.  Vector spaces of compact type}

\smallskip

In this section, we introduce a special  class of topological vector
spaces. This class arises naturally when considering locally analytic
representations, and it enjoys good properties with respect to passage
to subspaces, quotient spaces, and dual spaces.

We fix a spherically complete nonarchimedean field $K$ and let
$o_K$ denote its ring of integers. Define $\LC(K)$ to be the
category of locally convex $K$-vector spaces. We recall that a
topological $K$-vector space is called locally convex if it has a
fundamental system of open $0$-neighbourhoods consisting of
$o_K$-submodules. If $V$ and $W$ are objects in $\LC(K)$, we
denote by $\L(V,W)$ the space of continuous linear maps from $V$
to $W$. Following traditional usage, we write $\L_{s}(V,W)$ and
$\L_{b}(V,W)$ for the vector space $\L(V,W)$ equipped with its
weak and strong topologies respectively. For any $V$ in $\LC(K)$,
we denote by $V'$ the $K$-vector space $\L(V,K)$ and write
$V_{s}'$ for $\L_{s}(V,K)$ and $V_{b}'$ for $\L_{b}(V,K)$. Both of
these dual spaces are Hausdorff. An $o_K$-submodule $A$ of a
$K$-vector space $V$ is called a {\it lattice} if it $K$-linearly
generates $V$. Any open $o_K$-submodule of a locally convex
$K$-vector space is a lattice. If each closed lattice in $V$ is
open then the locally convex $K$-vector space $V$ is called {\it
barrelled}. This very basic property puts the Banach-Steinhaus
theorem into force ([B-TVS] II.25). We recall a few other concepts
before we introduce the notion in the title of this section.

\medskip

{\bf Definition:}

{\it Let $V$ and $W$ be Hausdorff locally convex $K$-vector spaces.}

{\it 1. A subset} $B \subseteq V$ {\it is called} $\ul {compactoid}$
{\it if for any open lattice} $A \subseteq V$ {\it there are finitely
many vectors} $v_1,\ldots,v_m \in V$ {\it such that} $B
\subseteq A + o_K v_1 + \ldots + o_Kv_m$.

{\it 2. An $o_K$-submodule} $B \subseteq V$ {\it is called}
$\ul{c-compact}$ {\it if it is compactoid and complete.}

{\it 3. A continuous linear map} $f:V\to W$ {\it is called} $\ul
{compact}$ {\it if there is an open lattice $A$ in $V$ such that
$f(A)$ is c-compact in $W$.}

{\it 4. $V$ is said to be} \ $\ul {of
\ compact \ type}$ {\it if it is the locally convex inductive limit of
a sequence}
$$
V_{1}{\buildrel\imath_{1}\over\lra}
V_{2}{\buildrel\imath_{2}\over\lra}V_{3}\lra\cdots
$$
{\it of locally convex Hausdorff $K$-vector spaces $V_{n}$, for}
$n\in\Ndss${\it , with injective compact linear maps} $\imath_{n}$.

\medskip

In part 3 of the above Definition one can, in fact, suppose that the
$V_n$ are Banach spaces (compare [GKPS] 3.1.4 or [Kom] Lemma 2).

\medskip

{\bf Theorem 1.1:}

{\it Any space $V$ of compact type is Hausdorff, complete,
bornological, and reflexive; its strong dual is a Fr\'{e}chet space
and satisfies $V_{b}'=\lim\limits_{\lla} (V_{n}')_{b}$.}

Proof: [GKPS] 3.1.7 or [Kom] Thm. 6' and Thm. 12 or [NFA] Lemma 7 and
Thm. 8 in III.5 (compare also [B-TVS] III.6 Prop.7).

\medskip

Let $\LC_{c}(K)$ be the full subcategory of $\LC(K)$ consisting of
spaces of compact type. We note that any reflexive space is
barrelled ([NFA] Lemma III.4.4).

\medskip

{\bf Proposition 1.2:}

{\it i. Let $V$ be a space of compact type, and let $U\subseteq V$ be
a closed vector subspace.  Then $U$ and $V/U$ are also of compact
type, and the maps in the exact sequence of vector spaces
$$
0\to (V/U)'_{b} \to V_{b}' \to U_{b}' \to 0
$$
are strict.

ii. The category $\LC_{c}(K)$ is closed under countable locally convex
direct sums (and therefore also finite products.)}

Proof: i. [Kom] Thm. 8 and Thm. 7', [GKPS] 3.1.16, and [B-TVS] IV.29
Cor. 1 and 2. ii. [Kom] Thm. 9 and 10.

\medskip

The duals of spaces of compact type have a natural characterization.

\medskip

{\bf Theorem 1.3:}

{\it A $K$-Fr\'{e}chet space $V$ is the strong dual of a space of
compact type if and only if $V$ is nuclear.}

Proof: [Kom] Thm. 17 or [GKPS] 2.5.8, 3.1.7, and 3.1.13.

\medskip

Let $\LC_{nF}(K)$ be the full subcategory of $\LC(K)$ consisting of
nuclear Fr\'{e}chet spaces.

\medskip

{\bf Corollary 1.4:}

{\it The functor
$$
\morphism{\LC_{c}(K)}{\LC_{nF}(K)}{V}{V_{b}'}
$$
is an anti-equivalence of categories.  For any two $V$ and $W$ in
$\LC_{c}(K)$, the natural linear map $\L_{b}(V,W)\to
\L_{b}(W_{b}',V_{b}')$ is a topological isomorphism.}

Proof: By reflexivity, applying the operation $V\mapsto V_{b}'$ twice
yields the identity functor.  The transpose operation gives a
continuous map from $\L_{b}(V,W)$ to $\L_{b}(W_{b}',V_{b}')$, which,
by reflexivity, has a continuous inverse.

\medskip

We conclude with a generalization of reflexivity for spaces of compact
type.  We will apply this result in our study of distributions in the
next section.

Suppose that  $V$ is of compact type, with defining sequence of Banach
spaces $V_{n}$, and that $W$ is a Banach space. An element $v\otimes
w$ of $V_{n}\otimes W$ defines an element of $\L((V_{n})'_b,W)$ by
$\ell \mapsto \ell(v)w$. $\L_b((V_{n})'_b,W)$ being again a Banach
space, we obtain a map $V_{n}\hat{\otimes}W\to\L_b((V_{n})'_b,W)$.

Composing with the natural map $V'_b=\lim\limits_{ \lla} (V_{n})'_b\to
(V_{n})'_b$ yields a map
$$
C:
\lim\limits_{\lra}
(V_{n} \hat{\otimes}  W)\to
\L(V'_b,W).
$$


{\bf Proposition 1.5:}

{\it The map $C$ is a linear isomorphism.}

Proof: For simplicity we denote the strong or Banach space dual of a
Banach space $U$ simply by $U'$. Because $V$ is of compact type, each
transition map $\imath_{n}:V_{n}\to V_{n+1}$ is compact. It follows
that the map $\imath_{n}'':V_{n}''\to V_{n+1}''$ factors through
$V_{n+1}\subseteq V_{n+1}''$ (see [Kom] Lemma 1 or [NFA] Lemma
III.5.3). Consequently, the inductive sequences
$$V_{1}\hat{\otimes} W\buildrel\imath_{1}\over\to  V_{2}\hat{\otimes} W
\buildrel\imath_{2}\over\to \cdots$$
and
$$V_{1}''\hat{\otimes} W\buildrel\imath_{1}''\over\to  V_{2}''\hat{\otimes} W
\buildrel\imath_{2}''\over\to \cdots\leqno{(1)}
$$
have the same limit.  Consider now the inductive sequence
$$
\L(V_{1}',W)\buildrel\jmath_{1}\over\to
\L(V_{2}',W)\buildrel\jmath_{2}\over\to \cdots\leqno{(2)}
$$
The transition maps  in this sequence are given by
$\jmath_{n}(f)=f\circ\imath_{n}'$.  Since the $\imath_{n}$ are
compact, so are the $\imath_{n}'$ (loc. cit.), and therefore the image
of $\L(V_{n}',W)$ under $\jmath_{n}$ lies in the subspace of compact
operators in $\L(V_{n+1}',W)$. By [Gru] Prop. 3.3.1, 5.2.2, and 5.3.3,
this space of compact operators is precisely $V_{n}''\hat{\otimes} W$.
As a result, the inductive sequences (1) and (2) have the same limit.
Now, because $W$ is a Banach space, any continuous linear map
$\lim\limits_{\lla} V_{n}'\to W$ must factor through one of the
$V_{n}'$. Therefore, we have a linear isomorphism
$$
\lim\limits_{\lra} \L(V_{n}',W)\buildrel\sim\over\to\L(\lim\limits_{\lla} V_{n}',
W)\ .
$$
By Thm. 1.1, $\lim\limits_{\lla} V_{n}'=V_{b}'.$  Tracing through the
steps in this isomorphism shows that the map is as described in the
above statement.

\bigskip

{\bf 2. Distributions}

\smallskip

In this section, we discuss the space of distributions on a
locally analytic mani-fold. We construct a general ``integration''
map allowing us to apply distributions to vector valued locally
analytic functions. We review the convolution product of
distributions on a locally analytic group. Finally, we show that
entire power series in Lie algebra elements converge to
distributions on the group.

We fix fields $\dQ\subseteq L\subseteq K$ such that $L/\dQ$ is finite
and $K$ is spherically complete with respect to a nonarchimedean
absolute value $|\; |$ extending the one on $L$. We let $o_L$ denote
the ring of integers in $L$. In the following $L$ plays the role of
the base field whereas $K$ will be the coefficient field.

Let $M$ be a $d$-dimensional locally $L$-analytic manifold
([B-VAR] \S 5.1 and [Fe1] \S 3.1). We always assume that $M$ is
strictly paracompact (and hence locally compact); this means that
any open covering of the topological space $M$ can be refined into
a covering by pairwise disjoint open subsets; e.g., if $M$ is
compact then it is strictly paracompact ([Fe1] \S 3.2).

Let $V$ be a Hausdorff locally convex $K$-vector space. In this
situation the locally convex $K$-vector space $C^{an} (M,V)$ of
all $V$-valued locally analytic functions on $M$ is defined ([Fe2]
2.1.10). We briefly recall this construction. A $V$-index $\cI$ on
$M$ is a family of triples
$$
\{(D_i,\phi_{i} ,V_i)\}_{i\in I}
$$
where the $D_i$ are pairwise disjoint open subsets of $M$ which
cover $M$, each $\phi_i :D_i\longrightarrow L^d$ is a chart of the
manifold $M$ whose image is an affinoid ball, and
$V_i\hookrightarrow V$ is a $BH$-space of $V$, i.e., an injective
continuous linear map from a $K$-Banach space $V_i$ into $V$. We
form the locally convex direct product
$$
\cF_\cI (V):=\prod\limits_{i\in I} \cF_{\phi_i} (V_i)
$$
of the $K$-Banach spaces
$$
\matrix{
\cF_{\phi_i} (V_i):= & \hbox{all functions}\;\; f:D_i\longrightarrow V_i\;\;
 \hbox{such that}\;\; f\circ\phi^{-1}_{i}\;\; \hbox{is a}\hfill\cr
 &  V_i\hbox{-valued holomorphic function on the affinoid ball}\cr
 & \phi_i (D_i).\hfill\cr}
$$
Denoting by $\cO(\phi_{i} (D_i))$ the $K$-Banach algebra of all
power series in $d$ variables with coefficients in $K$ converging
on all points of the affinoid $\phi_{i} (D_i)$ defined over an
algebraic closure of $K$ we have
$$
\cF_{\phi_i} (V_i)\cong\cO (\phi_i (D_i)) \mathop{\hat{\otimes}}\limits_K V_i \ .
$$
The $V$-indices on $M$ form a directed set on which $\cF_\cI (V)$
is a direct system of locally convex $K$-vector spaces. We then
may form the locally convex inductive limit
$$
C^{an}
(M,V):=\mathop{\lim}\limits_{\mathop{\longrightarrow}\limits_\cI}
\ \cF_\cI (V) \ .
$$
The formation of this space is compatible with disjoint coverings
in the following sense ([Fe2] 2.2.4): Whenever
$M=\mathop{\bigcup}\limits^{\cdot}_{j\in J} M_j$ is an open
covering by pairwise disjoint subsets $M_j$ then one has the
natural topological isomorphism
$$
C^{an} (M,V) =\prod\limits_{j\in J} C^{an} (M_j ,V) \ .
$$

\medskip

{\bf Definition:}

{\it The strong dual $D(M,K):=C^{an} (M,K)_b'$ of $C^{an} (M,K)$
is called the locally convex vector space of $K$-valued
distributions on $M$.}

\medskip

{\bf Lemma 2.1:}

{\it If $M$ is compact then $C^{an} (M,K)$ is of compact type; in
particular, $D(M,K)$ is the projective limit over the $K$-indices
on $M$ of the strong dual spaces $\cF_\cI (K)_b'$ and is a
$K$-Fr\'echet space.}

Proof: Apply [ST] Lemma 1.5 (in fact, a straightforward
generalization to the case of a spherically complete $K$) to the
transition maps in the direct system $\cF_\cI (K)$ or use [Fe2]
2.3.2.

\medskip

For a general $M$ we may choose a covering
$M=\mathop{\bigcup}\limits^{\cdot}_{j\in J} M_j$ by pairwise
disjoint compact open subsets and obtain a topological isomorphism
$$
D(M,K)=\mathop{\oplus}\limits_{j\in J} D(M_j,K) \ .
$$
This observation together with the above Lemma shows that the
strong topology on $D(M,K)$ coincides with the topology considered
in [Fe1] 2.5.5 and 3.4.1.

The only completely obvious elements in $D(M,K)$ are the Dirac
distributions $\delta_x$, for $x \in M$, defined by $\delta_x
(f):=f(x)$.

\medskip

{\bf Theorem 2.2:} (Integration)

{\it If $V$ is the union of countably many BH-spaces then the map
$$
\matrix{
I^{-1}:\L (D(M,K),V) & \buildrel\cong\over\longrightarrow & C^{an}
(M,V)\cr\cr
\hfill A & \longmapsto & [x\longmapsto A(\delta_x )]\cr}
$$
is a well defined $K$-linear isomorphism.}

Proof: Clearly the map in the assertion is compatible with
disjoint open coverings of $M$. We therefore may assume that $M$
is compact so that $D(M,K)$ is a Fr\'echet space. On the other
hand, since the sum of two $BH$-spaces again is a $BH$-space we
find, by our assumption on $V$, an increasing sequence
$V_1\subseteq V_2\subseteq\ldots$ of $BH$-spaces of $V$ such that
$V=\mathop{\bigcup}\limits_{n\in \N} V_n$. By the closed graph
theorem ([B-TVS] I. 20 Prop. 1) any continuous linear map from the
Fr\'echet space $D(M,K)$ into $V$ factors through some $V_n$. In
other words we have
$$
\L (D(M,K),V)=\mathop{\lim}\limits_{\mathop{\longrightarrow}\limits_{n\in\N}}
\ \L (D(M,K),V_n) \ .
$$
Moreover, by the same closed graph theorem any $BH$-space of $V$
is contained in some $V_n$. Since $M$ is compact (so that in the
definition of $C^{an} (M,V)$ we need to consider only $V$-indices
whose underlying covering of $M$ is finite) this means that
$$
C^{an}
(M,V)=\mathop{\lim}\limits_{\mathop{\longrightarrow}\limits_{n\in\N}}
\ C^{an} (M,V_n) \ .
$$
Hence we are further reduced to the case that $V$ is a Banach
space. But then we are precisely in the situation of Prop. 1.5.
The only additional observation to make is that the map
$$
\matrix{
I:\cO (\phi_i (D_i))\mathop{\hat{\otimes}}\limits_K V &
\longrightarrow & \L (\cO (\phi_i (D_i))^{'}_{b},V)\cr\cr
\hfill f\otimes \upsilon & \longmapsto & [\lambda\longmapsto\lambda (f)\upsilon ]\hfill\cr}
$$
considered there satisfies $I(f\otimes\upsilon ) (\delta_x )
=f(x)\cdot\upsilon$ for $x\in M$.

\medskip

Since any locally analytic function on a compact manifold factors
through a $BH$-space the arguments in the above proof show that
for an arbitrary $V$ we still have a natural map
$$
I:C^{an}(M,V) \longrightarrow \L (D(M,K),V)
$$
such that $I(f) (\delta_x )=f(x)$ for any $f\in C^{an} (M,K)$ and
$x\in M$. This map $I$ should be considered as integration: Given
a locally analytic function $f:M\longrightarrow V$ and a
distribution $\lambda$ on $M$ we may formally write
$$
I(f)(\lambda )=\int\limits_M f(x) d\lambda (x).
$$

\smallskip

The case where $M:=G$ is a locally $L$-analytic group has
additional features which we now want to discuss. First of all we
note that any such $G$ is strictly paracompact ([Fe1] 3.2.7).
According to [Fe1] 4.1.6 there is a unique separately continuous
$K$-bilinear map
$$
\hat{\otimes} :D(G,K)\times D(G,K)\longrightarrow D(G\times G,K)
$$
with the property
$$
(\lambda \hat{\otimes}\mu) (f_1\times f_2)=\lambda (f_1)\cdot\mu
(f_2)
$$
for any $\lambda ,\mu\in D(G,K),f_i\in C^{an} (G,K)$, and
$(f_1\times f_2) (g_1,g_2):=f_1 (g_1)\cdot f_2 (g_2)$. If $G$ is
compact then $\hat{\otimes}$ is (jointly) continuous by [Fe1]
4.2.1 (the separate continuity in general is an immediate
consequence of that). The group multiplication $m:G\times
G\longrightarrow G$ induces by functoriality a continuous linear
map
$$
m_\ast :D(G\times G,K)\longrightarrow D(G,K).
$$
We define the convolution $\ast$ on $D(G,K)$ by
$$
\ast := m_\ast\circ\hat{\otimes} \ ,
$$
i.e.,
$$
(\lambda\ast\mu) (f)=(\lambda\hat{\otimes}\mu) (f\circ m)\;\;
\hbox{for}\;\; \lambda ,\mu\in D(G,K)\;\; \hbox{and}\;\; f\in
C^{an} (G,K).
$$

\medskip

{\bf Proposition 2.3:}

{\it $(D(G,K),\ast)$ is an associative $K$-algebra with the Dirac
distribution $\delta_1$ in $1\in G$ as the unit element; the
convolution $\ast$ is separately continuous. If $G$ is compact
then $D(G,K)$ is a Fr\'echet algebra.}

Proof: [Fe1] 4.4.1 and 4.4.4.

\medskip

The selfmap $g\longmapsto g^{-1}$ of $G$ induces an
anti-automorphism of $D(G,K)$. For this reason we never have to
consider right modules for $D(G,K)$ and we will simply speak of
modules when we mean left modules.

\smallskip

One method to explicitly construct elements in $D(G,K)$ is through
the Lie algebra $\gfr$ of $G$. Let $C^{\omega}_{1}$ be the
``stalk'' of $C^{an}(G,K)$ at the identity element, as defined in
[Fe2] \S2.3.1.  This is the locally convex inductive limit of the
spaces $C^{an}(Y,K)$ as $Y$ runs through the family of compact
open neighborhoods of the identity in $G$. Note that each map
$C^{an}(Y,K)\to C^{\omega}_{1}$ is surjective.

For each compact open subgroup $H \subseteq G$ the Lie algebra
$\gfr$ of $G$ acts on $C^{an}(H,K)$  by continuous endomorphisms
defined by
$$
  (\xfr f)(g):={d\over dt} f(\exp (-t\xfr ) g)_{|t=0}
$$
or written in a more invariant way, by using the left translation
action of the group $H$ on $C^{an}(H,K)$,
$$
\xfr f := \mathop{\lim}\limits_{t\to 0} {{\rm exp}(t\xfr)f-f \over t}
$$
for $\xfr\in\gfr$ where $\exp:\gfr\, {--->}\, G$ denotes the
exponential map defined locally around 0 ([Fe2] 3.1.2 and 3.3.4). This
extends, by the universal property, to a left action of the universal
enveloping algebra $U(\gfr)$ on each $C^{an}(H,K)$ (and hence on
$C^{\omega}_{1}$) by continuous endomorphisms.  In particular, an
element $\zfr$ of $U(\gfr)$ gives a continuous linear form on any
$C^{an}(H,K)$ (and hence on $C_{1}^{\omega}$) by the rule
$$
f\mapsto (\dot{\zfr}(f))(1)\;\;.
$$
where $\zfr \mapsto \dot{\zfr}$ is the unique anti-automorphism of
$U(\gfr)$ which extends the multiplication by $-1$ on $\gfr$ ([B-GAL]
Chap.I,\S2.4). By [Fe2] 4.7.4 this induces an injection $U(\gfr)
\to (C_{1}^{\omega})'$. Moreover, the continuous surjection
$C^{an}(G,K)\to C_{1}^{\omega}$ yields a continuous injection
$(C_{1}^{\omega})'_b\to D(G,K)$. Altogether we have the chain of
inclusions
$$
U(\gfr) \mathop{\otimes}\limits_{L} K \hookrightarrow
(C_{1}^{\omega})'
\hookrightarrow D(G,K) \;\;.
$$
The space $C_{1}^{\omega}$ can be represented explicitly. Fix an
ordered basis $\xfr_{1},\ldots,\xfr_{d}$ of $\gfr$. The corresponding
``canonical system of coordinates of the second kind'' ([B-GAL]
Chap.III,\S4.3) is the map
$$
\morphism{{\rm small}\;{\rm ball}\;{\rm around}\; 0 \; {\rm in}\; L^d}
{G}{(t_{1},\ldots,t_{d})}
     {\exp(t_{1}\xfr_{1})\cdots\exp(t_{d}\xfr_{d})\;\;.}
$$
This chart gives us, in a neighborhood of $1$, and for a
multi-index $\ul{n}$, the monomial functions
$$
T^{\ul{n}}(\exp(t_{1}\xfr_{1})\cdots\exp(t_{d}\xfr_{d})):=t_{1}^{n_{1}}\cdots
t_{d}^{n_{d}}\;\;.
$$

\medskip

{\bf Lemma 2.4:}

{\it The space} $C_{1}^{\omega}$ {\it is isomorphic to the space
of all power series}
$$
C_{1}^{\omega}=
\{\ \sum_{\ul{n}} a_{\ul{n}}T^{\ul{n}}\ : a_{\ul{n}}\in K
\hbox{\rm\ and\ }|a_{\ul{n}}|r^{|\ul{n}|}\to 0
 \hbox{\rm\ as $|\ul{n}|\to\infty$ for some $r>0$} \};
$$
{\it its dual is the space of all power series}
$$
(C_{1}^{\omega})'=\{\ \sum_{\ul{n}}b_{\ul{n}}Z^{\ul{n}}\ :
b_{\ul{n}}\in K
\hbox{\rm\ and $\mathop{\sup}\limits_{\ul{n}}|b_{\ul{n}}|r^{-|\ul{n}|} < \infty$ for any
$r>0$}\};
$$
{\it the strong topology on} $(C_{1}^{\omega})'_b$ {\it is defined
by the family of norms}
$$
\|\sum_{\ul{n}}
b_{\ul{n}}Z^{\ul{n}}\|_{r}:=\mathop{\sup}\limits_{\ul{n}}
|b_{\ul{n}}|r^{-|\ul{n}|}\;\;\;{\rm for}\;\; r > 0 \;;
$$
{\it finally, the map} $U(\gfr)\to (C_{1}^{\omega})'$ {\it is
given explicitly by}
$$
\xfr_{1}^{n_{1}}\cdots \xfr_{d}^{n_{d}}\mapsto
(-1)^{|\ul{n}|}\ul{n}!Z^{\ul{n}}\;\;.
$$

Proof:  By definition, $C_{1}^{\omega}$ consists of functions
represented by  power series with a non-zero radius of convergence
around $1$. This is precisely the condition we specify here.  The
description of the dual follows easily, with $Z^{\ul{n}}$ dual to
$T^{\ul{n}}$.  The map on $U(\gfr)$ comes from the proof of Lemma
4.7.2 in [Fe2] which shows that
$$
(((\xfr_{1}^{i_{1}}\cdots
\xfr_{d}^{i_{d}})^{\bf{\cdot}})T^{\ul{n}})(1)=(-1)^{|\ul{n}|}\ul{n}!
\delta_{\ul{i},\ul{n}}
$$
where
$$
|\ul{n}| = \sum_{j=1}^{d} n_{j} \hbox{\rm\ \ and\ \ }
\ul{n}! = n_{1}!n_{2}!\cdots n_{d}!
$$
and $\delta_{\ul{i},\ul{n}}$ is Kronecker's symbol.

\medskip

{\bf Lemma 2.5:}

{\it Let} $\{b_{\ul{n}}\}_{\ul{n}}$ {\it be a set of elements of}
$K$ {\it such that} $\mathop{\sup}\limits_{\ul{n}}
|b_{\ul{n}}|r^{-|\ul{n}|} < \infty$ {\it for any} $r > 0${\it ;
then the series}
$$
\sum_{\ul{n}} b_{\ul{n}}\xfr_{1}^{n_{1}}\cdots \xfr_{d}^{n_{d}}
$$
{\it converges in} $(C_{1}^{\omega})'_b$ {\it (that is, its
sequence of partial sums converges to an element of}
$(C_{1}^{\omega})'_b$ {\it ).}

Proof: Given the previous Lemma, it suffices to show that the
series $\sum _{\ul{n}}b_{\ul{n}}Z^{\ul{n}}$ is the limit of its
partial sums.  We compute:
$$
\eqalign{
\|\sum _{\ul{n}}b_{\ul{n}}Z^{\ul{n}}-\sum_{|\ul{n}|\leq n}
b_{\ul{n}}Z^{\ul{n}}\|_{r}&=
\mathop{\sup}\limits_{|\ul{n}|>n} |b_{\ul{n}}|r^{-|\ul{n}|}\cr
&=\mathop{\sup}\limits_{|\ul{n}|>n}
|b_{\ul{n}}|\left({{r}\over{2}}\right)^{-|\ul{n}|}\left({{1}\over{2}}\right)^{|\ul{n}|}\cr
&\le \left({{1}\over{2}}\right)^{n}
\|\sum_{\ul{n}}b_{\ul{n}}Z^{\ul{n}}\|_{r/2}\cr &\to 0\hbox{\rm\ as
$n\to\infty$} }
$$

\medskip

{\bf Corollary 2.6:}

{\it Given any set of coefficients} $\{b_{\ul{n}}\}_{\ul{n}}$ {\it
as in the Lemma, the series}
$$
\sum_{\ul{n}} b_{\ul{n}}\xfr_{1}^{n_{1}}\cdots \xfr_{d}^{n_{d}}
$$
{\it converges in} $D(G,K)$.

Proof: This follows from the continuity of the map
$(C_{1}^{\omega})'_b\to D(G,K)$.

\medskip

In somewhat vague terms this Corollary says that any entire power
series over $K$ in the $\xfr_1,\ldots,\xfr_d$ converges in
$D(G,K)$.

It is useful to have various other interpretations of the above
inclusion of $\gfr$ into $D(G,K)$. The $\gfr$-action on $C^{an}(G,K)$
induces a $\gfr$-action on $D(G,K)$ by $(\xfr \lambda)(f) := \lambda
((-\xfr)f)$. In particular $\xfr \delta_1$ is the image of $\xfr$ in
$D(G,K)$. One easily computes
$$
(\xfr \lambda)(f) = \mathop{\lim}\limits_{t\to 0}  {{\rm
exp}(t\xfr)\lambda - \lambda \over t}(f) \ \ .
$$
This means that in $C^{an}(G,K)'_s$ we have the limit formula
$$
\xfr \lambda = \mathop{\lim}\limits_{t\to 0}\  ({\delta_{{\rm
exp}(t\xfr)} - \delta_1 \over t}*\lambda) \ \ .
$$
As a direct product of spaces of compact type the space
$C^{an}(G,K)$ is reflexive. For such a space $C^{an}(G,K)'_s$
coincides with $C^{an}(G,K)'_b$ given its weakened topology, and
any bounded subset in $C^{an}(G,K)'_s$ (e.g., a convergent
sequence in this space together with its limit) is contained in a
bounded and c-compact $o_K$-submodule of $C^{an}(G,K)'_b = D(G,K)$
([Ti1] III \S3 and \S4 and [Ti2] Thm. 2 or [NFA] Satz 2 in III.3
and Satz 3 and Lemma 4 in III.4). Moreover, on submodules of the
latter type the strong topology coincides with the weak topology
([Gru] 5.3.4 or [NFA] Satz 8 in III.1). It follows that the above
limit formula even holds in $D(G,K)$. But the convolution in
$D(G,K)$ is separately continuous. We therefore obtain
$$
\xfr \lambda = (\mathop{\lim}\limits_{t\to 0} {\delta_{{\rm
exp}(t\xfr)} - \delta_1 \over t})*\lambda = (\xfr \delta_1)*\lambda \
\ {\rm in} \ D(G,K) \ .
$$

\bigskip

{\bf 3. Locally analytic $G$-representations}

\smallskip

In this section, we show that locally analytic $G$-representations on
$K$-vector spaces $V$ of compact type  are equivalent to a certain
class of modules for the ring $D(G,K)$.  We let $\dQ
\subseteq L
\subseteq K$ be fields as in the previous section and we let $G$
be a $d$-dimensional locally $L$-analytic group with Lie algebra
$\gfr$.

\medskip

{\bf Definition:}

{\it A locally analytic $G$-representation $V$
(over $K$) is a barrelled locally convex Hausdorff $K$-vector
space $V$ equipped with a $G$-action by continuous linear
endomorphisms such that, for each $v\in V$, the orbit map
$\rho_{v}(g):=gv$ is a $V$-valued locally analytic function on
$G$.}

\medskip

By the Banach-Steinhaus Theorem, the map $(g,v)\to gv$ is
continuous, as is the  Lie algebra action
$$
v\to {\xfr}v:={{d}\over{dt}}\exp(t\xfr)v|_{t=0}
$$
for $\xfr\in \gfr$  and $v\in V$.  This $\gfr$-action extends to
an action of the universal enveloping algebra $U(\gfr)$ on $V$ by
continuous linear endomorphisms. Taylor's formula (compare [Fe2]
3.1.4) says that, for each fixed $v\in V$ there is a  a
sufficiently small neighborhood $U$ of 0 in $\gfr$ such that, for
$\xfr\in U$, we have a convergent expansion
$$
\exp(\xfr)v=\sum_{n=0}^{\infty} {{1}\over{n!}}\xfr^{n}v\ .
$$

Now suppose that $V$ is a locally analytic representation of $G$.
We obtain a candidate for a $D(G,K)$-module structure on $V$ by
defining
$$
\morphism{D(G,K)\times V}{V}{(\mu,v)}{\mu*v=I(\rho_{v})(\mu)}
$$
where $I$ is the integration map from the previous section. For a
Dirac distribution $\delta_{g}$, we have $\delta_{g}*v=gv$.

\medskip

{\bf Lemma 3.1:}

{\it The Dirac distributions generate a dense subspace in
$D(G,K)$.}

Proof: If $H \subseteq G$ is a compact open subgroup then we have
$$
   D(G,K) = \mathop{\bigoplus}\limits_{g \in G/H} \delta_g *
   D(H,K)\ \ .
$$
This observation shows that we may assume $G$ to be compact. Let
$\Delta$ be the closure in $D(G,K)$ of the subspace generated by
the Dirac distributions. Let $\ell$ be a linear form on $D(G,K)$
vanishing on $\Delta$. By reflexivity of $D(G,K)$, such a linear
form must be given by a locally analytic function $f$ on $G$, and
to say $\ell$ vanishes on $\Delta$ is to say that $f$ is
identically zero on $G$.  By the Hahn-Banach Theorem, we conclude
that $\Delta =D(G)$.

\medskip

{\bf Proposition 3.2:}

{\it The map $(\mu,v)\to\mu*v$ is separately continuous in $\mu$
and $v$, and $V$ is a module over $D(G,K)$; this $D(G,K)$-module
structure extends the action of $U(\gfr)$ on $V$. Any continuous
linear $G$-map between locally analytic $G$-representations is a
$D(G,K)$-module homomorphism.}

Proof: The continuity in $\mu$ is clear by construction. To show
the continuity in $v$ we use the previous Lemma. If $G$ is compact
then $D(G,K)$ is metrisable and any $\mu\in D(G,K)$ is
approximated by a sequence $\mu =\lim_{n\to\infty} \mu_{n}$ where
each $\mu_{n}$ is a linear combination of Dirac distributions.
This remains true for an arbitrary $G$ by the same observation as
in the previous proof. Now, each $\mu_{n}$ is a continuous linear
map from $V$ to $V$, and $\mu$ is their pointwise limit. Since $V$
barrelled, we see from the Banach-Steinhaus theorem that $\mu$
itself is a continuous linear map from $V$ to $V$.  To show that
$V$ is in fact a $D(G,K)$-module, we must show that, for fixed
$v$, $\mu'*(\mu*v)=(\mu'*\mu)*v$. This identity clearly holds for
Dirac distributions and it follows for general $\mu$ by
continuity. The argument for the last assertion is similar.
Finally, since the action of $\gfr$ and $U(\gfr)$ is expressible
as a limit of actions by Dirac distributions, the remaining claim,
too, follows from continuity.

\medskip

If $V$ is a barrelled Hausdorff space which is the union of
countably many BH-spaces then the above Proposition has a
converse: As a consequence of Theorem 2.2 any separately
continuous $D(G,K)$-module structure on $V$ comes from a locally
analytic $G$-representation on $V$ in the way described above. In
other words, the category of locally analytic $G$-representations
on such spaces is equivalent to the category of separately
continuous $D(G,K)$-module structures on them.

\medskip

{\bf Corollary 3.3:}

{\it The functor}
$$
\matrix{
\matrix{
\hbox{\rm locally analytic $G$-representations}\cr
\hbox{\rm on $K$-vector spaces of compact type}\cr
\hbox{\rm with continuous linear $G$-maps}\cr
\cr
} &\to&
\matrix{
\hbox{\rm separately continuous $D(G,K)$-}\cr
\hbox{\rm modules on nuclear Fr\'echet}\cr
\hbox{\rm spaces with continuous}\cr
\hbox{\rm $D(G,K)$-module maps}\cr
}\cr &&\cr V&\mapsto&V'_{b}\cr}
$$
{\it is an anti-equivalence of categories.}

Proof: This follows from the above discussion together with Cor. 1.4
once we show that $V$ carries a separately continuous $D(G,K)$-module
structure if and only $V'_b$ does. Indeed, the separate continuity of
the pairing $D(G,K)\times V\to V$ is equivalent to the assertion that
the map $D(G,K)\to\L_s(V,V)$ is continuous. It follows from ([B-TVS]
III.31 Prop.6) that this map remains continuous when $\L_b(V,V)$ is
given its strong topology. Applying Cor. 1.4, we obtain a continuous
map from $D(G,K)$ into $\L_{b}(V_{b}',V_{b}')$, so that $V_{b}'$
carries a separately continuous $D(G,K)$-module structure. By the
reflexivity of $V$ this argument works the same way in the opposite
direction.

\medskip

If the group $G$ is compact then $D(G,K)$ is a Fr\'echet space. Since
any separately continuous bilinear map between Fr\'echet spaces is
jointly continuous ([B-TVS] III.30 Cor. 1) we may reformulate in this
case the last assertion as follows.

\eject

{\bf Corollary 3.4:}

{\it If $G$ is compact then the functor}
$$
\matrix{
\matrix{
\hbox{\rm locally analytic $G$-representations}\cr
\hbox{\rm on $K$-vector spaces of compact type}\cr
\hbox{\rm with continuous linear $G$-maps}\cr
} &\to&
\matrix{
\hbox{\rm continuous $D(G,K)$-modules on}\cr
\hbox{\rm nuclear Fr\'echet spaces with con-}\cr
\hbox{\rm tinuous $D(G,K)$-module maps}\cr
}\cr &&\cr V&\mapsto&V'_{b}\cr}
$$
{\it is an anti-equivalence of categories.}

\medskip

Let $\M_K(G)$ be the category of all left $D(G,K)$-modules in the
purely algebraic sense. The full subcategory of $\M_K(G)$
generated by the objects which are countably generated (i.e., such
that there exists a surjective $D(G,K)$-module map
$\oplus_{\Ndss}D(G,K) \to M$) will be denoted by $\M_{K}^{c}(G)$.

\medskip

{\bf Definition:}

{\it Assume $G$ to be compact; a locally analytic
$G$-representation $V$ is called admissible if $V$ is of compact
type and $V'_b$ as a $D(G,K)$-module is countably generated.}

\medskip

For a compact group $G$ we let $Rep^{adm}_K(G)$ denote the
category of all admissible $G$-representations with continuous
linear $G$-maps.

\medskip

{\bf Lemma 3.5:}

{\it Assume $G$ to be compact; the category $Rep^{adm}_K(G)$ is closed
with respect to the passage to closed $G$-invariant subspaces.}

Proof: Let $U \subseteq V$ be a closed $G$-invariant subspace. By
Prop. 1.2.i, $U$ is of compact type. Using [Fe2] 1.2.3 and the fact
that the coefficients of the power series expansions of the orbit maps
$\rho_v$ are limits of linear combinations of values it easily follows
that $U$ is a locally analytic $G$-representation. Finally, looking at
the exact sequence in Prop. 1.2.i we see that $U'$ is a quotient
module of the countably generated $D(G,K)$-module $V'$ and hence is
countably generated as well.

\medskip

Later on we will make use of the following simple consequence of this
Lemma: Any $V$ in $Rep^{adm}_K(G)$ which is such that $V'$ is an
(algebraically) simple $D(G,K)$-module must be topologically
irreducible as a $G$-representation.

\eject

{\bf Proposition 3.6:}

{\it If $G$ is compact then the contravariant functor
$$
\matrix{
Rep^{adm}_K(G)  &\to& \M^c_K(G) \cr &&\cr V&\mapsto&V'\cr}
$$
is a fully faithful embedding.}

Proof: Let $M$ and $N$ be $D(G,K)$-modules in the image of the
functor (in fact, $N$ can be any separately continuous
$D(G,K)$-module). We have to show that any abstract
$D(G,K)$-module homomorphism $\alpha : M
\to N$ is continuous. Choose generators $(m_{i})_{i\in\Ndss}$ for
$M$, and let
$$
\beta :\oplus_{\Ndss}D(G,K)\to M
$$
be the map $\beta ((\mu_{i}))= \sum_{\Ndss} \mu_{i}m_{i}$. By
([B-TVS] II.34, Cor. to Prop.10), this map is strict; in other
words, $M$ carries the quotient topology with respect to the map
$\beta$. At the same time, we have a map
$$
\gamma :\oplus_{\Ndss}D(G,K)\to N
$$
given by $\gamma ((\mu_{i}))=\sum_{\Ndss} \mu_{i}\alpha(m_{i})$.
Since the $D(G,K)$-action on $N$ is continuous, this map is
continuous, and by construction it factors through $M$, where it
induces the original map $\alpha$.  By the universal property of
the quotient topology, $\alpha$ is continuous.

\medskip

This result can be seen as a first step towards algebraizing the
theory of admissible $G$-representations.

A further important algebraic concept is concerned with the action
of the centre $Z(\gfr)$ of the universal enveloping algebra
$U(\gfr)$.

\medskip

{\bf Proposition 3.7.:}

{\it If $G$ is an open subgroup of the group of $L$-rational
points of a connected algebraic $L$-group $\bf G$ then $Z(\gfr)
\mathop{\otimes}\limits_{L} K$ lies in the centre of $D(G,K)$.}

Proof: Since $\gfr = {\rm Lie}(G) = {\rm Lie}(\bf G)$ and $U(\gfr)
\subseteq D(G,K) \subseteq D({\bf G}(L),K)$ we may assume that $G
= {\bf G}(L)$. By Lemma 3.1 it suffices to show that $\delta_g *
\zfr * \delta_{g^{-1}} = \zfr$ for any $g \in {\bf G}(L)$ and any
$\zfr \in Z(\gfr)$. This amounts to ([B-GAL] III \S3.12 Cor. 2)
the adjoint action of ${\bf G}(L)$ on $U(\gfr)$ fixing the centre
$Z(\gfr)$ which is a consequence of the connectedness of $\bf G$
([DG] II \S6.1.5).

\medskip

Let us assume for the rest of this section that $G$ satisfies the
assumption of the above Proposition.

\medskip

{\bf Definition:}

{\it 1. A $D(G,K)$-module $M$ is called $Z(\gfr)$-locally finite
if $Z(\gfr)m$, for each $m \in M$, is finite dimensional over $K$.

2. A locally analytic $G$-representation $V$ is called
$Z(\gfr)$-locally cofinite if the $D(G,K)$-module $V'$ is
$Z(\gfr)$-locally finite.}

\medskip

Clearly, a finitely generated $D(G,K)$-module $M$ which is
$Z(\gfr)$-locally finite actually is $Z(\gfr)$-finite in the sense
that the annihilator ideal
$$
\afr_M := \{\zfr \in Z(\gfr)\mathop{\otimes}\limits_L K : \zfr M =
0\}
$$
is of finite codimension in $Z(\gfr)\mathop{\otimes}\limits_L K$.
If $M$ is a $Z(\gfr)$-locally finite and simple $D(G,K)$-module
then $\afr_M$ is a maximal ideal in
$Z(\gfr)\mathop{\otimes}\limits_L K$: By Schur's lemma ${\rm
End}_{D(G,K)}(M)$ is a skew field which implies that
$Z(\gfr)\mathop{\otimes}\limits_L K/\afr_M$ is an integral domain.
But being of finite dimension $Z(\gfr)\mathop{\otimes}\limits_L
K/\afr_M$ then has to be a field. In this case the action of
$Z(\gfr)$ on $M$ therefore is given by a character $\chi_M$ of
$Z(\gfr)$ into the multiplicative group of the finite field
extension $Z(\gfr)\mathop{\otimes}\limits_L K/\afr_M$ of $K$. This
$\chi_M$ is called the $\underline{infinitesimal\ character}$ of
$M$. If $K$ is algebraically closed then the infinitesimal
character is a character $\chi_M : Z(\gfr) \longrightarrow
K^{\times}$. Assume that $K$ is algebraically closed, $G$ is
compact, and $M$ is a simple $D(G,K)$-module which is the dual $M
= V'$ of an admissible and $Z(\gfr)$-locally finite
$G$-representation $V$ then $Z(\gfr)$ acts on $V$ through the
character $\dot{\chi}_M : Z(\gfr)
\longrightarrow K^{\times}$ defined by $\dot{\chi}_M (\zfr) :=
\chi_M (\dot{\zfr})$.

\bigskip

{\bf 4. The one dimensional case}

\smallskip

As the first important example we will discuss in this section the
case of the group $G = \dZ$. Hence in the following $L = \dQ$, and
$G$ always denotes the one dimensional additive group $\dZ$. The
spherically complete field $K/\dQ$ is assumed to be contained in
$\dC$.

In order to determine the structure of the ring $D(\dZ,K)$ we use the
theory of the Fourier transform. The character group
$$
\hat{G}:=\Hom_{an} (G,\Cdss^{\times}_p)
$$
of $G= \dZ$ is defined to be the group of all locally $\dC$-analytic
group homomorphisms $\kappa :G\longrightarrow \Cdss_p^{\times}$. If
$X$ denotes the rigid analytic open unit disk over $K$ then it is easy
to see that
$$
\matrix{
\hat{G} & \longrightarrow & X(\dC)\hfill\cr
\kappa & \longmapsto & \kappa (1)-1\cr}
$$
is a well defined injective map. It is in fact bijective since, for
any $z\in X(\dC)$ and any $a\in\dZ$ the series
$$
\kappa_z (a):= (1+z)^a:=\sum\limits_{n\ge 0} {a\choose n} z^n
$$
converges and defines a locally $\dC$-analytic character $\kappa_z$ of
$\dZ$ with $\kappa_z (1)-1=z$.  Actually $\kappa_z$, for $z\in X(K)$,
is locally $K$-analytic so that we have the embedding
$$
\matrix{
X(K) & \hookrightarrow & C^{an}(\dZ,K)\cr
\hfill z & \mapsto & \kappa_z \hfill \cr}\ .
$$
Hence any linear form $\lambda\in D(\dZ,K)$ gives rise to the function
$$
F_\lambda (z):=\lambda (\kappa_z)
$$
on $X(K)$ which is called the \underbar{Fourier transform} of
$\lambda$. In fact, $F_\lambda$ is a holomorphic function on $X$. Let
$\cO (X)$ denote the ring of all $K$-holomorphic functions on $X$. We
recall that
$$
\matrix{\cO (X) =& \hbox{all power series}\ \ F(T) =\sum\limits_{n\ge 0} a_n
T^n\ \ \hbox{with}\ a_n\in K\cr & \hbox{which converge on}\ X(\dC)\ .
\hfill}
$$
For any $r \in p^{\Qdss}$ let $X_r$ denote the $K$-affinoid disk of
radius $r$ around zero. The ring $\cO (X_r)$ of all $K$-holomorphic
functions on $X_r$ is a $K$-Banach algebra with respect to the
multiplicative spectral norm. Since
$$
\cO (X)=\mathop{\lim}\limits_{\mathop{\longleftarrow}\limits_{r<1}}\
 \cO (X_r) =\mathop{\bigcap}\limits_{r<1} \cO (X_r)
$$
the ring $\cO (X)$ is an integral domain and is, in a natural way, a
$K$-Fr\'echet algebra.

\medskip

{\bf Theorem 4.1:} (Amice)

{\it The Fourier transform
$$
\matrix{
D(\dZ,K)& \buildrel\cong\over\longrightarrow & \cO (X)\cr
\hfill \lambda & \longmapsto & F_\lambda\cr}
$$
is an isomorphism of $K$-Fr\'echet algebras.}

Proof: [Am2] 1.3 and 2.3.4 (based on [Am1]); compare [Sch] for a
concise write-up of this proof.

\medskip

By work of Lazard the structure of the ring $\cO (X)$ is well known.
We recall that an effective divisor on $X_r$ is a map
$D:X_r\longrightarrow \Z_{\ge 0}$ with finite support. The effective
divisors are partially ordered by
$$
D\le D'\;\; \hbox{iff}\;\; D(x)\le D'(x)\;\; \hbox{for any}\;\; x\in
X_r\ .
$$
Moreover, for any $F \in \cO (X_r)$ there is the associated divisor of
zeros $(F)$ on $X_r$.

We define
$$
\matrix{
{\rm Div}^+ (X) := & \hbox{all maps}\;\; D:X\longrightarrow\Zdss_{\ge
0}\;\; \hbox{such that, for any}\;\; r<1,\cr &
\hbox{the restriction}\;\; D|X_r\;\; \hbox{has finite support}\hfill\cr}
$$
to be the partially ordered set of effective divisors on $X$ and
we obtain the divisor map
$$
\matrix{
\cO (X)\setminus \{ 0\} & \longrightarrow & {\rm Div}^+ (X)\cr
\hfill F & \longmapsto & (F) \hfill \cr}
$$
as the projective limit of the divisor maps for the $X_r$ with $r<1$.
It is clear that, for any $0\not= F,F'\in\cO (X)$, we have
$$
F|F'\;\; \hbox{if and only if}\;\; (F)\le (F')\ ;
$$
in particular:
$$
(F)=(F')\;\; \hbox{if and only if}\;\; F=uF'\;\; \hbox{for some}\;\;
u\in\cO (X)^{\times}\ .
$$

\medskip

{\bf Theorem 4.2:} (Lazard)

{\it i. The divisor map $\cO (X)\setminus
\{ 0\}\longrightarrow{\rm Div}^+ (X)$ is surjective;

ii. the map
$$
\matrix{
{\rm Div}^+ (X) & \buildrel\sim\over\longrightarrow & \hbox{\rm {all
nonzero closed ideals in}}\; \cO (X)\cr\cr
\hfill D & \longmapsto & I_D:=\{f\in\cO (X):(f)\ge D\}\cup\{ 0\}\cr}
$$
is a well defined bijection;

iii. the closure of a nonzero ideal $I \subseteq \cO (X)$ is the ideal
$I_{D(I)}$ for the divisor $D(I)(x):=\mathop{\min}\limits_{0\not= f\in
I}\ (f)(x)$;

iv. in $\cO (X)$ the three families of closed ideals, finitely
generated ideals, and principal ideals coincide.}

Proof: [Laz] Thm. 7.2, Prop. 7.10, (7.3), and Prop. 8.11 (again we
refer to [Sch] for a concise write-up).

\medskip

We record the following facts for later use.

\medskip

{\bf Lemma 4.3:}

{\it The translation action $\lambda\mapsto {^b\lambda}$, for $b\in
\dZ$, of $\dZ$ on $D(\dZ,K)$ corresponds under the Fourier transform to
$$
F_{({^b\lambda})}(T)=(1+T)^b F_{\lambda}(T) \ .
$$
If $\xfr_1 \in {\rm Lie}(\dZ) = \dQ$ corresponds to $1 \in \dQ$ then
the Lie algebra action $\lambda \mapsto \xfr_1\lambda$ on $D(\dZ,K)$
corresponds to
$$
F_{\xfr_1 \lambda}(T)=\log(1+T)F_{\lambda}(T) \ .
$$
For $b \in \dZ$, the multiplication by $b$ on $\dZ$ induces an
operator $m_{b*}$ on $D(\dZ,K)$ which corresponds to
$$
F_{m_{b*}\lambda}(T)=F_{\lambda}((1+T)^{b}-1) \ .
$$
Finally, the operator $\Theta(f) := af(a)$ on $C^{an}(\dZ,K)$ induces
the map
$$
F_{\lambda\circ\Theta}(T)=(1+T){{dF_{\lambda}(T)}\over{dT}} \ .
$$
}

Proof: These are easy computations (see also [Am2] 2.3).

\medskip

Having determined the structure of the ring $D(\dZ,K)$ we now turn to
a discussion of the module theory for this ring. More specifically we
are aiming at some information about the image of our fully faithful
embedding
$$
\matrix{
Rep^{adm}_K(\dZ)  &\to& \M^c_K(\dZ) \cr\cr  V&\mapsto&V'\cr}
$$
in Prop. 3.6.

\medskip

{\bf Proposition 4.4:}

{\it Let $M$ be a free $D(\dZ,K)$-module of finite rank equipped with
the obvious direct product topology; then any finitely generated
$D(\dZ,K)$-submodule $N$ of $M$ is closed.}

Proof: The ring $D(\dZ,K) = \cO (X)$ is an adequate ring in the sense
of [Hel] \S2 since:\hfill\break
 - it is an integral domain;\hfill\break
 - any finitely generated ideal is principal (by Thm.
 4.2.iv);\hfill\break
 - for any two elements $F \neq 0$ and $F'$ in $\cO (X)$ the
 ''relatively prime part $R$ of $F$ w.r.t. $F'\ $'' exists: Define the
 divisor $D$ by
 $$
 D(x) := \left\{\matrix{(F)(x) & \hbox{if}\ (F')(x)=0\;,\cr
 \ 0\hfill & \hbox{otherwise}\hfill\cr}\right.
 $$
 and let $R$ be such that $(R) = D$ (by Thm. 4.2.i).\hfill\break
  The main result Thm. 3 in [Hel] then says that the elementary divisor theorem holds
 over $\cO (X)$. This means that there is an isomorphism $M\cong \cO
 (X)^m$ under which $N$ corresponds to a submodule of the form
 $F_1\cO (X)\oplus\ldots\oplus F_n\cO (X)$ with $n\leq m$ and $F_i
 \in\cO (X)$. Since, by Thm. 4.2.iv, principal ideals are closed in $\cO
 (X)$ it follows that $N$ is closed in $M$.

\medskip

Let now $M$ be any finitely presented $D(\dZ,K)$-module. If we choose
a finite presentation
$$
D(\dZ,K)^n \rightarrow D(\dZ,K)^m \rightarrow M \rightarrow 0
$$
then as a consequence of the Proposition the quotient topology on $M$
is Hausdorff. Hence $M$ with this topology is a nuclear Fr\'echet
space with a continuous $D(\dZ,K)$-action. Moreover, this quotient
topology is independent of the chosen presentation. We therefore see
that any $D(\dZ,K)$-module of finite presentation $M$ carries a
natural Fr\'echet topology such that $M \cong V'_b$ for some
admissible $\dZ$-representation $V$.

By Lazard's theorem the simple $D(\dZ,K)$-modules of finite
presentation are (up to isomorphism) those of the form $D(\dZ,K)/\mfr$
for some closed maximal ideal $\mfr \subseteq D(\dZ,K)$ and they all
are finite dimensional as $K$-vector spaces. Any such is the dual of
an irreducible admissible $\dZ$-representation.

\medskip

{\bf Remark:}

{\it A topologically irreducible locally analytic
$\dZ$-representation $V$ is finite dimensional (and irreducible)
if and only if it is $Z(\dQ)$-locally cofinite.}

Proof: The other implication being trivial we assume that $V$ is
$Z(\dQ)$-locally cofinite. Fix a nonzero $\ell \in V'$. Then the
closed ideal $I := \{F \in \cO (X) : F*\ell = 0\}$ in $\cO (X) =
D(\dZ,K)$ is nonzero and proper. Choose a point $x_0 \in X$ such
that $D(I)(x_0) > 0$. The closed ideal $J \subseteq \cO (X)$ such
that $D(J)(x) = D(I)(x)$ for $x \neq x_0$ and $D(J)(x_0) =
D(I)(x_0)-1$ has the property that $J/I$ is finite dimensional.
Hence $\{v
\in V : \ell '(v) = 0\ \hbox{for any}\ \ell ' \in J*\ell\}$ is a
$D(\dZ,K)$-invariant proper closed subspace of finite codimension
in $V$. By the topological irreducibility of $V$ such a subspace
must be zero. Hence $V$ is finite dimensional.

\bigskip

{\bf 5. The principal series of the Iwahori subgroup in} $GL_2(\dQ)$

\smallskip

In this section we let $G = B$ be the Iwahori subgroup of $GL_2(\dQ)$.
This is the subgroup of $GL_2(\dZ)$ consisting of all matrices which
are lower triangular mod $p$. As before, $K/\dQ$ is a spherically
complete extension field contained in $\dC$. Let $P_{\rm o}^-$, resp.
$T_{\rm o}$, denote the subgroup of $B$ of all upper triangular, resp.
diagonal, matrices.

We fix a $K$-valued locally analytic character
$$
\chi:T_{\rm o} \to K^{\times}
$$
and we define $c(\chi)\in K$ so that
$$
\chi(\pmatrix{t^{-1} & 0\cr 0 & t})=\exp(c(\chi)\log(t))
$$
for $t$ sufficiently close to $1$ in $\dZ$.

Our aim is to study the locally analytically induced
$B$-representation
$$
\Ind_{P_{\rm o}^-}^{B}(\chi):=\{ f\in C^{an}(B,K): f(gp)=
\chi(p^{-1})f(g) \ \hbox{for any}\ g\in B,
 p\in P^{-}_{\rm o}\}
$$

where $B$ acts by left translation.

\medskip

{\bf Proposition 5.1:}

{\it $\Ind_{P_{\rm o}^-}^{B}(\chi)$ is a locally analytic
$B$-representation; the underlying vector space is topologically
isomorphic to $C^{an}(\dZ,K)$ and in particular is of compact type.}

{\bf Proof:}   That $\Ind_{P_{\rm o}^-}^{B}(\chi)$  is locally
analytic follows from [Fe2] 4.1.5 and the fact that $B/P_{\rm o}^-$ is
compact. The map
$$
\matrix{
\morphism{\iota:\dZ}{B}{b}{\pmatrix{1 & 0\cr b & 1\cr}}\cr}
$$
is a section of the projection map $B\to B/P_{\rm o}^-=\dZ$, and
applying [Fe2] 4.3.1, we see that pullback by this section yields the
asserted isomorphism of topological vector spaces. The last claim then
follows from Lemma 2.1.

\medskip

{\bf Definition:}

{\it Let $M^{-}_{\chi}:=\Ind_{P_{\rm o}^-}^{B}(\chi)'_{b}$ be the
$D(B,K)$-module obtained from $\Ind_{P_{\rm o}^- }^{B}(\chi)$ by
applying the duality functor of Cor. 3.4.}

\medskip

As a vector space, the module $M^{-}_{\chi}$ is isomorphic to
$D(\dZ,K)$, and the  $D(\dZ,K)$-module structure on $M^{-}_{\chi}$
coming from the ring inclusion
$$
\iota_{*}:D(\dZ,K)\to D(B,K)
$$
is simply the ring multiplication.  Thus $M^{-}_{\chi}$ is a copy of
$D(\dZ,K)$ with the additional structure corresponding to the action
of $D(B,K)$. Before we calculate this additional structure explicitly
we remark that this observation already implies the admissibility of
the $B$-representation $\Ind_{P_{\rm o}^- }^{B}(\chi)$.

\medskip

{\bf Lemma 5.2:}

{\it Under the identification $\Ind_{P_{\rm
o}^-}^{B}(\chi)=C^{an}(\dZ,K)$ we have:

i. The translation action of $\pmatrix{1&0\cr 0 &b\cr} \in B$ is given
by the formula
$$
\pmatrix{1&0\cr 0 &b\cr}f = \chi(\pmatrix{1&0\cr 0
&b\cr})\cdot(m_{b^{-1}}^* f) \ ;
$$
ii. the translation action of $\pmatrix{1&b\cr 0 &1\cr} \in B$ is
given by
$$
(\pmatrix{1 & b\cr 0 & 1\cr}f)(a)=\chi(\pmatrix{(1-ab)^{-1} & 0\cr 0
&(1-ab)\cr}) f({{a}\over{1-ab}}) \ ;
$$
iii. the Lie algebra element $\ufr^{+}:=\pmatrix{0 & 0\cr 1&0\cr}\in
\gfr\lfr_2(\dQ) = \Lie(B)$ acts by
$$
(\ufr^{+}f)(a)= -{{df(a)}\over{da}} \ ;
$$
iv. the Lie algebra element $\ufr^{-}:=\pmatrix{0 & 1\cr 0&0\cr}\in
\Lie(B)$ acts by
$$
\ufr^{-}f = -c(\chi)\Theta(f) - \Theta^2(\ufr^+ f) \ ;
$$
here $m_{b^{-1}}^*$ and $\Theta$ are the operators from Lemma 4.3.}

Proof:  An elementary computation.

\medskip

{\bf Lemma 5.3:}

{\it Let
$$
c_{i}^{(n)}={{n!}\over{i!}}\pmatrix{{c(\chi)-i}\cr{n-i}}
\ \ for \ 0 \leq i \leq n \ ;
$$
then, after making the identification $M^{-}_{\chi} = D(\dZ,K)$ and
applying the Fourier transform, the action of $(\ufr^-)^m \in U({\rm
Lie}(B))$, for each $m \geq 0$, satisfies
$$
F_{(\ufr^{-})^{m}\lambda} = \sum_{i=0}^{m} (-1)^{i} c_{i}^{(m)}{\rm
(log}(1+T))^i\Delta^{m+i}F_{\lambda}
$$
where the operator $\Delta$ is defined by $\Delta F :=
(1+T){{dF(T)}\over{dT}}$. Moreover, the translation action of
$\pmatrix{1 & 0\cr 0&b\cr}\in B$ satisfies}
$$
F_{({1\atop 0}{0\atop b})\lambda}(T) = \chi (\pmatrix{1 & 0\cr
0&b^{-1}\cr})\cdot F_{\lambda}((1+T)^b-1) \ .
$$

Proof: Using the formulas
$$c_i^{(m)}(c(\chi) - m - i) + c_{i-1}^{(m)} = c_i^{(m+1)}
$$
and
$$
\ufr^+ \Theta^m = -m\Theta^{m-1} + \Theta^m \ufr^+
$$
as well as Lemma 5.2.iv one obtains by induction the identity
$$
(\ufr^{-})^{m} = (-1)^m \sum_{i=0}^{m}
c_{i}^{(m)}\Theta^{m+i}\circ(\ufr^{+})^i
$$
in $\L(C^{an}(\dZ,K),C^{an}(\dZ,K))$. By Lemma 4.3 this transposes
into the first of the asserted identities. The second one is a direct
consequence of Lemma 5.2 and Lemma 4.3.

\medskip

We now come to the main result of this section.

\medskip

{\bf Theorem 5.4:}

{\it If $c(\chi)\not\in\Ndss_{0}$, then $M^{-}_{\chi}$ is an
(algebraically) simple $D(B,K)$-module.}

Proof: As we pointed out above, $M^{-}_{\chi}$ is isomorphic to
$D(\dZ,K)$, and, under this isomorphism, the subalgebra $D(\dZ,K)
\subseteq D(B,K)$ acts via the usual multiplication action of $D(\dZ,K))$ on
itself. Let $I$ be a nonzero $D(B,K)$-submodule of $M^{-}_{\chi}$;
then $I$ corresponds to a nonzero $D(B,K)$-invariant ideal in
$D(\dZ,K)$. Thus we may apply Lazard's theory of divisors to study
$I$. (Note, however, that we do NOT assume that $I$ is closed.)

In a first step we claim that the ideal $I$ is generated by elements
$F$ such that $(F)$ is supported on the set
$$
\mu_{\infty}:=\{z \in \dC : 1+z\hbox{\rm\ is a root
of unity}\}.
$$
To see this claim, choose a nonzero $F\in I$, let $S$ be the set of
zeros of $F$ and let $S_{0}:=S\backslash S\cap\mu_{\infty}$. For any
pair $z,z'\in S_{0}$, there is at most one $u\in\Zdss^{\times}_p$ such
that
$$
u_*(z) := (1+z)^{u}-1=z'\ .
$$
Because $\Zdss^{\times}_p$ is uncountable, and $S_{0}$ is countable,
there must be a $u\in\Zdss^{\times}_p$ so that $u_{*}(S_{0})\cap
S_{0}=\emptyset$. As a result of the second part of the previous
Lemma, the common support of the divisors $(F)$ and $(\pmatrix{1 &
0\cr 0 & u}F)$ is contained in $\mu_{\infty}$. By Thm. 4.2 the ideal
$<F,\pmatrix{1 & 0\cr 0 & u}F>$ is generated by an element whose
divisor is supported on $\mu_{\infty}$.  Since $F\in I$ was arbitrary,
and $I$ is translation invariant, $I$ must be generated by elements
whose divisors are supported on $\mu_{\infty}$.

If, as before $X$ denotes the open unit disk considered as a rigid
analytic variety over $K$ we let the sequence of points
$x_0,x_1,\ldots$ in $X$ correspond to the orbits of the absolute
Galois group of $K$ in the set $\mu_{\infty}$.

We now begin with a nonzero element $F\in I$ with the following
property:

{
\parindent=20pt

\item{(*)} There is a sequence of positive integers $0<m_{0}<m_{1}<\ldots$
such that the divisor $(F)$ is supported on $x_0,x_1,\ldots$ and, for
each integer $k\geq 0$, the multiplicity of $(F)$ in $x_k$ is $m_k$.

}

The existence of such an $F$ follows from the above first step, and
from the principality of all divisors (see Thm. 4.2.i). We will now
use the action of $\ufr^{+}$ to perturb the element $F$, thereby
obtaining a new element $F^c \in I$ whose divisor is supported
entirely off the set $\mu_{\infty}$. By Lazard's theory, $F$ and $F^c$
together generate $D(\dZ,K)$, and we will conclude that $I=D(\dZ,K)$.
This perturbation will be carried out by induction. Choose a function
$\sigma:\Ndss_0\to p^{\Qdss}$ such that, for all real numbers $r>0$,
$\sigma(n)r^{n}\to 0$ as $n\to\infty$. First, observe that by the
first part of the previous Lemma
$$
(\ufr^{+})^{m_{k}}F\equiv [\prod_{i=0}^{m_k-1}(c(\chi)-i)]\cdot
\left((1+T){{d}\over{dT}}\right)^{m_{k}}F\pmod{\log(1+T)\cO (X)}\ .
$$
Because the divisor $(\log(1+T))$ has multiplicity one in each $x_k$,
we see that the divisor of $(\ufr^{+})^{m_{k}}F$ has multiplicity zero
in $x_k$ but has positive multiplicity in each $x_j$ for $j>k$
(because of $m_{j}>m_{k}$). Now suppose that we have found elements
$b_{0},\ldots, b_{k-1}\in K$ so that $|b_{i}|<\sigma(i)$ for $0\leq i
< k$ and the divisor of
$$
F_{k-1}:=(b_{0}(\ufr^{+})^{m_{0}}+\cdots+b_{k-1}(\ufr^{+})^{m_{k-1}})F
$$
has zero multiplicity in $x_0,\ldots,x_{k-1}$. Note, however, that
$(F_{k-1})$ must have positive multiplicity in $x_k$, while
$(\ufr^{+})^{m_{k}}F$ does not. By avoiding finitely many special
cases, we may choose $b_{k}\in K$ so that $|b_{k}|<\sigma(k)$ and the
divisor of $F_{k-1}+b_{k}(\ufr^{+})^{m_{k}}F$ has zero multiplicity in
$x_0,\ldots,x_k$. Having constructed the $b_k$ in this way for all
$k\geq 0$ our choice of $\sigma$, together with Cor. 2.6, tells us
that the element
$$
\ufr_{\infty}:=\sum_{k=0}^{\infty} b_{k}(\ufr^+)^{m_{k}}\in D(B,K)
$$
is well defined so that we have $F^c:=\ufr_{\infty}F \in I$. Since, by
construction, the divisor of
$(b_{0}(\ufr^{+})^{m_{0}}+\cdots+b_{k}(\ufr^{+})^{m_{k}})F$ has zero
multiplicity in $x_k$ whereas $\sum_{j>k} b_{j}(\ufr^+)^{m_{j}}F$ has
positive multiplicity it follows that $(F^c)$ is supported completely
outside of $\mu_{\infty}$. Hence we see that $I=D(\dZ,K)$. This proves
that $M^{-}_{\chi}$ is simple.

\medskip

In order to formulate the next result we introduce the character
$\epsilon$ defined by $\epsilon(\pmatrix{t_{1} & 0\cr 0 & t_{2}}):=
t_{2}/t_{1}$.

\medskip

{\bf Proposition 5.5:}

{\it For any two $K$-valued locally analytic characters $\chi \neq
\chi '$ of $\ T_{\rm o}$ the vector space ${\rm
Hom}_{D(B,K)}(M^{-}_{\chi '},M^{-}_{\chi})$ of all $D(B,K)$-module
homomorphisms is $1$-dimensional if $c(\chi) \in \Ndss_{0}$ and $\chi
' =
\epsilon^{-1-c(\chi)} \chi$ and is zero otherwise.}

Proof:  By Prop. 3.6 we may dualize back and determine the
$B$-equivariant continuous linear maps from $\Ind^{B}_{P^{-}_{\rm
o}}(\chi)$ to $\Ind^{B}_{P^{-}_{\rm o}}(\chi ')$. By Frobenius
reciprocity ([Fe2] 4.2.6) this is equivalent to determining the space
of all $P^{-}_{\rm o}$-equivariant continuous linear maps from
$\Ind^{B}_{P^{-}_{\rm o}}(\chi)$ to $K_{\chi '}$ where the latter
denotes the $1$-dimensional representation of $P^{-}_{\rm o}$ given by
$\chi '$. Making our usual identification $M^{-}_{\chi} \cong
D(\dZ,K)$ this identifies with the space
$$
 \Jscr := \{F \in D(\dZ,K): \chi '(g)\cdot gF = F\ \hbox{\rm for any}\ g
 \in P^{-}_{\rm o}\}\ .
$$
We study this defining condition separately for three different kinds
of elements $g$ in $P^{-}_{\rm o}$.

1. Any $g = \pmatrix{b & 0\cr 0 & b}$ acts on $\Ind^{B}_{P^{-}_{\rm
o}}(\chi)$ by multiplication with $\chi (g)$ and hence on $D(\dZ,K)$
by multiplication with $\chi^{-1} (g)$. For $\Jscr$ to be nonzero we
therefore must have
$$
\chi \mid \Zdss^{\times}_p \cdot \pmatrix{1 & 0\cr 0 & 1} = \chi '
\mid \Zdss^{\times}_p \cdot \pmatrix{1 & 0\cr 0 & 1}\ . \leqno(1)
$$
2. According to Lemma 5.3 any $g = \pmatrix{1 & 0\cr 0 & b}$ acts on
$D(\dZ,K)$ via
$$
(gF)(T) = \chi(\pmatrix{1 & 0\cr 0 & b^{-1}})\cdot F((1+T)^b -1)\ .
$$
Any nonzero $F \in \Jscr$ therefore satisfies
$$
F(T) = {\chi ' \over \chi}(\pmatrix{1 & 0\cr 0 & b}) \cdot F((1+T)^b
-1)\ . \leqno(2)
$$
Applying the operator ${d\over db}|_{b=1}$ to this equation we obtain
the differential equation
$$
{c(\chi ')-c(\chi) \over 2}F(T)+(1+T){\rm log}(1+T){dF(T) \over dT} =
0\ . \leqno(3)
$$
Here we have used that
$$
2{d \over db}({\chi ' \over \chi}(\pmatrix{1 & 0\cr 0 & b}))|_{b=1} =
{d \over db}({\chi ' \over \chi}(\pmatrix{b^{-1} & 0\cr 0 &
b}))|_{b=1} = c(\chi ')-c(\chi)\ .
$$
The differential equation $(3)$ has a nonzero solution in $D(\dZ,K)$
if and only if
$$
c(\chi)-c(\chi ') \in 2\Ndss_{0}\ . \leqno(4)
$$
In order to see this put $c := (c(\chi )-c(\chi '))/2$ and assume that
$F(T) =\sum\limits_{n\ge 0} a_n T^n$ is a solution of $(3)$. By
inserting $F$ into $(3)$ and comparing coefficients we obtain
$$
ca_0 = 0\ \ \hbox{and}\ \ (n-c)a_n = \sum\limits_{i=1}^{n-1} (-1)^i
{n-i \over i(i+1)} a_{n-i} \ \ \hbox{for}\ \ n \geq 1\ .
$$
From this one easily deduces by induction:\hfill\break
 - in case $c \not\in \Ndss_0$ that $a_n = 0$ for all $n \geq 0$,
 i.e., that $F = 0$;\hfill\break
 - in case $ c \in \Ndss_0$ that $a_n = 0$ for $n<c$ and that $a_n$ for
 $n>c$ can be recursively computed from $a_c$, i.e., that up to
 scalars there is at most one $F$.\hfill\break
 In the latter case one verifies that
$$
F(T) := [{\rm log}(1+T)]^{(c(\chi)-c(\chi '))/2}
$$
indeed is a solution of $(3)$. One immediately checks that this
solution satisfies the equation $(2)$ if and only if
$$
\chi '(\pmatrix{1 & 0\cr 0 & b}) = b^{(c(\chi ')-c(\chi))/2}\cdot
\chi(\pmatrix{1 & 0\cr 0 & b})\ \ \ \hbox{for any}\ b \in
\Zdss^{\times}_p\ .\leqno(5)
$$
Note that the conditions $(1),(4)$ and $(5)$ together are equivalent
to the requirement that
$$
\chi ' = \epsilon^{-m}\chi \ \ \ \hbox{for some}\ m \in \Ndss_{0}\
. \leqno(6)
$$
3. Assuming that this latter identity holds we finally look at the
action of $g = \pmatrix{1 & b\cr 0 & 1}$. We again consider first the
derived action of $\ufr^{-}$. By Lemma 5.3 the solution $F(T) := ({\rm
log}(1+T))^m$ must satisfy the equation
$$
c(\chi)\Delta F = {\rm log}(1+T)\Delta^2 F
$$
if it lies in $\Jscr$. But this means (since $\chi ' \neq \chi$) that
$c(\chi) = m-1 \geq 0$ which together with $(6)$ is precisely the
condition in our assertion. It remains to show that under this
condition $({\rm log}(1+T))^m$ indeed is invariant under any $g =
\pmatrix{1 & b\cr 0 & 1}$. This we do by looking at the distribution
$$
\morphism{C^{an}(\dZ,K)}{K}{f(a)}{(({d \over da})^m f)(0)}
$$
which by the second assertion in Lemma 4.3 and Lemma 5.2.iii
corresponds to $({\rm log}(1+T))^m$ under the Fourier transform. By
induction we deduce from Lemma 5.2.ii the formula
$$
(({d \over da})^m (gf))(a) = \sum_{i=0}^m (-1)^{m+i} c_i^{(m)}
{b^{m-i} \over (1-ab)^{m+i}}(g(({d \over da})^i f))(a)
$$
with the $c_i^{(m)}$ as defined in Lemma 5.3. But since $c(\chi)= m-1$
we actually have $c_i^{(m)} = 0$ for $i<m$. Setting $a=0$ and using
Lemma 5.2.ii again we therefore obtain
$$
(({d \over da})^m (gf))(0) = (({d \over da})^m f)(0)\ .
$$

\medskip

Letting $P^{+}_{\rm o} \subseteq B$ denote the subgroup of lower
triangular matrices we may analogously consider the $B$-representation
by left translation on the induction $\Ind^{B}_{P^{+}_{\rm o}}(\chi)$
from $P^{+}_{\rm o}$. It is admissible and, setting $M^{+}_{\chi} :=
\Ind^{B}_{P^{+}_{\rm o}}(\chi)'_{b}$, we have the following result
whose proof is completely parallel to the proofs of Thm. 5.4 and Prop.
5.5 and is therefore omitted.

\medskip

{\bf Theorem 5.6:}

{\it i. If $c(\chi) \not\in -\Ndss_{0}$ then $M^{+}_{\chi}$ is an
(algebraically) simple $D(B,K)$-module;

ii. for any two $\chi \neq \chi '$ the vector space ${\rm
Hom}_{D(B,K)}(M^{+}_{\chi '},M^{+}_{\chi})$ is $1$-dimensional if
$c(\chi) \in -\Ndss_{0}$ and $\chi ' =
\epsilon^{1-c(\chi)} \chi$ and is zero otherwise.}

\eject

{\bf Corollary 5.7:}

{\it If $c(\chi) \not\in \pm\Ndss_{0}$ then $\Ind^{B}_{P^{\mp}_{\rm o}
}(\chi)$ is topologically irreducible as a $B$-representation.}

\medskip

{\bf Proposition 5.8:}
$$
{\rm Hom}_{D(B,K)}(M^{+}_{\chi},M^{-}_{\chi '}) = {\rm
Hom}_{D(B,K)}(M^{-}_{\chi},M^{+}_{\chi '}) = 0 \ \ \ \hbox{\it for any
$\chi$\ \hbox{and} $\chi '$}\ .
$$
Proof: We only discuss the vanishing of the first space, the other
case being completely analogous. By the same reasoning as at the
beginning of the proof of Prop. 5.5 we have to show that there are
no nontrivial $P^{+}_{\rm o}$-equivariant continuous linear maps
from $\Ind^{B}_{P^{-}_{\rm o}}(\chi ')$ to $K_{\chi}$. Each such
is a continuous linear form on $\Ind^{B}_{P^{-}_{\rm o}}(\chi ')$
which is invariant under the unipotent subgroup $\iota (\dZ)$ of
$P^{+}_{\rm o}$. Using Prop. 5.1 it is sufficient to show that
there is no nonzero distribution in $D(\dZ,K)$ which is
$\dZ$-invariant. But this is the well known nonexistence of a
''$p$-adic Haar measure'' on $\dZ$ (one way to see this is to use
the first assertion in Lemma 4.3).

\medskip

Define the elements $\efr := \pmatrix{1 & 0\cr 0 & 1}$ and $\hfr
:= \pmatrix{1 & 0\cr 0 & -1}$ in $\gfr\lfr_2 (\dQ)$ and
$$
\cfr := {1 \over 2}\hfr^2 + \ufr^+ \ufr^- + \ufr^- \ufr^+ \in
U(\gfr\lfr_2 (\dQ))\ .
$$
It is well known that the centre of $U(\gfr\lfr_2 (\dQ))$ is the
polynomial ring
$$
Z(\gfr\lfr_2 (\dQ)) = \dQ [\efr ,\cfr ]
$$
in the two variables $\efr, \cfr$. Using the formulas iii. and iv.
in Lemma 5.2 and the corresponding formula
$$
(\hfr f)(a) = -c(\chi)f(a) + 2a{df(a) \over da}
$$ one easily checks that $\cfr$ acts on $M_{\chi}^{\pm}$ through
multiplication by the scalar
$$
({1 \over 2}c(\chi) \mp 1)c(\chi) \ .
$$
In particular, $M_{\chi}^{\pm}$ is $Z(\gfr\lfr_2 (\dQ))$-finite.

\bigskip

{\bf 6. The principal series of} $GL_2(\dZ)$ {\bf and} $GL_2(\dQ)$

\smallskip

The results of the preceding section will enable us in this section to
analyze completely the irreducibility properties of the principal
series of the group $GL_{2}(\dQ)$.

We let $G:=GL_2(\dQ)$ and $G_{\rm o}:=GL_2(\dZ)$. Furthermore, $P$
denotes the Borel subgroup of lower triangular matrices in $G$, $T$
the subgroup of diagonal matrices, and $B$ as before the Iwahori
subgroup in $G_{\rm o}$. The field $K$ is as before.

This time we fix a $K$-valued locally analytic character
$$
\chi:T\to K^{\times}
$$
and we again define $c(\chi)\in K$ so that
$$
\chi(\pmatrix{t^{-1} & 0\cr 0 & t})=\exp(c(\chi)\log(t))
$$
for $t$ sufficiently close to $1$ in $\dZ$. The corresponding
principal series representation is $\Ind_{P}^{G}(\chi)$ with $G$
acting, as always, by left translation. From the Iwasawa decomposition
$G=G_{\rm o}P$ and [Fe2] 4.1.4 we obtain that the obvious restriction
map
$$
\Ind_{P}^{G}(\chi) \mathop{\longrightarrow}\limits^{\simeq}
\Ind_{P\cap G_{\rm o}}^{G_{\rm o}}(\chi)
$$
is a $G_{\rm o}$-equivariant topological isomorphism.(For simplicity
we denote by $\chi$ also the restriction of the original $\chi$ to
various subgroups of $T$.) Using the Bruhat decomposition
$$
G_{\rm o} = B(G_{\rm o}\cap P) \mathop{\cup}\limits^{\cdot} Bw(G_{\rm
o}\cap P)
$$
with $w:= \pmatrix{0 & -1\cr 1 & 0} \in G_{\rm o}$ and [Fe2] 2.2.4 and
4.1.4 we may further decompose into
$$
\Ind_{P\cap G_{\rm o}}^{G_{\rm o}}(\chi) = \Ind^{B}_{P^{+}_{\rm
o}}(\chi) \oplus \Ind^{B}_{P^{-}_{\rm o}}(w\chi) \ .
$$
This is a $B$-equivariant decomposition of topological vector spaces.
First of all we see that, as a consequence of the previous section,
$\Ind_{P\cap G_{\rm o}}^{G_{\rm o}}(\chi)$ is an admissible $G_{\rm
o}$-representation.

\medskip

{\bf Definition:}

{\it Let $M_{\chi}:=\Ind_{P}^{G}(\chi)'_{b}$ be the $D(G,K)$-module
obtained from $\Ind_{P}^{G}(\chi)$ by applying the duality functor of
Cor. 3.3.}

\medskip

{\bf Theorem 6.1:}

{\it If $c(\chi)\not\in -\Ndss_{0}$, then $M_{\chi}$ is an
(algebraically) simple $D(G_{\rm o},K)$-module. In parti-cular,
$\Ind^{G}_{P}(\chi)$ is a topologically irreducible $G_{\rm o}$- and a
fortiori $G$-representation.}

Proof: As a $D(B,K)$-module $M_{\chi}$ is isomorphic to the direct sum
$$
M_{\chi} = M_{\chi}^{+} \oplus M_{w\chi}^{-}
$$
of the two modules $M_{\chi}^{+}$ and $M_{w\chi}^{-}$ which by Thm.
5.4, Thm. 5.6.i, and Prop. 5.8 are simple and nonisomorphic under our
assumption on $c(\chi)$. Hence these two are the only nonzero proper
$D(B,K)$-submodules of $M_{\chi}$. Since the action of $w$ mixes the
two summands it follows that $M_{\chi}$ has to be simple as a
$D(G_{\rm o},K)$-module.

\medskip

{\bf Proposition 6.2:}

{\it  For any two $\chi \neq \chi '$ we have
$$
{\rm Hom}_{D(G,K)}(M_{\chi '},M_{\chi}) = {\rm Hom}_{D(G_{\rm
o},K)}(M_{\chi '},M_{\chi})
$$
and this vector space is $1$-dimensional if $c(\chi) \in
-\Ndss_{0}$ and $\chi ' = \epsilon^{1-c(\chi)} \chi$ and is zero
otherwise.}

Proof: It is immediate from 5.5, 5.6, and 5.8 that ${\rm
Hom}_{D(G_{\rm o},K)}(M_{\chi '},M_{\chi})$ is at most
$1$-dimensional and is zero unless $c(\chi) \in -\Ndss_{0}$ and
$\chi ' = \epsilon^{1-c(\chi)} \chi$. Let us henceforth assume
that $m := c(\chi) \in -\Ndss_{0}$ and that $\chi ' =
\epsilon^{1-m} \chi$. In order to establish our assertion it
remains to exhibit a nonzero $G$-equivariant continuous linear map
$$
I : \Ind_{P}^{G}(\chi) \longrightarrow \Ind_{P}^{G}(\chi').
$$
We observe that $U(\gfr \lfr_2 (\dQ))$ acts continuously on
$C^{an}(G,K)$ by left invariant differential operators which,
solely for the purposes of this proof, we denote by $(\zfr , f)
\mapsto \zfr f$; for $\xfr \in \gfr \lfr_2 (\dQ)$ the formula is
$$
(\xfr f)(g):={d\over dt} f(g\exp (t\xfr ))_{|t=0}\ .
$$
We claim that the operator on $C^{an}(G,K)$ corresponding to the
element $\zfr := (\ufr^- )^{1-m}$ restricts to a linear map $I$ as
above, i.e., that the map
$$
\matrix{I : \Ind_{P}^{G}(\chi) & \longrightarrow &
\Ind_{P}^{G}(\chi')\cr\cr
\hfill f & \longmapsto & (\ufr^-)^{1-m} f\hfill}
$$
is well defined. For this we have to show that, given an $f \in
\Ind_{P}^{G}(\chi)$, we have
$$
(\zfr f)(gp)=\chi'(p^{-1})\cdot (\zfr f)(g)\ \ \hbox{for\ any}\ g
\in G, p \in P.
$$
We check this separately for diagonal matrices and for lower
triangular unipotent matrices.  First let $h$ be a diagonal
matrix. Then
$$
(\zfr f)(gh)=
\chi(h^{-1})\epsilon^{1-m}(h^{-1})\cdot (\zfr f)(g)
=\chi'(h^{-1})\cdot (\zfr f)(g)
$$
as required, using the easily checked fact that
$\Ad(h)(\ufr^-)=\epsilon(h^{-1})\ufr^-$. Next we consider a lower
triangular unipotent matrix
$$
u=\pmatrix{1 & 0\cr a & 1}=\exp(a\ufr^{+})\ .
$$
Since $f$ is right $u$-invariant we have
$$
(\zfr f)(gu) = ((\Ad(u)\zfr)f)(g)\ .
$$
Now
$$
\eqalign{
\Ad(u)(\zfr) & =  \Ad(\exp(a\ufr^{+}))(\zfr)=\exp(\ad(a\ufr^{+}))(\zfr)\cr
& =  \zfr+a\ad(\ufr^+)(\zfr) +
(a^2/2)\ad(\ufr^+)(\ad(\ufr^+)(\zfr)) + \ldots}
$$
where the ``series'' in this expression is actually finite. Using
the elementary identity
$$
\ad(\ufr^{+})(\zfr) = (m-1)(\ufr^{-})^{-m}(-m+\hfr)
$$
in $U(\sfr\lfr_2 (\dQ))$ we see that
$$
[\ad(\ufr^+)]^j (\zfr) \in U(\gfr\lfr_2 (\dQ))(-m+\hfr ) +
U(\gfr\lfr_2 (\dQ))\ufr^+
$$
for any $j \geq 1$. Hence
$$
\Ad(u)(\zfr) \in \zfr + U(\gfr\lfr_2 (\dQ))(-m+\hfr ) +
U(\gfr\lfr_2 (\dQ))\ufr^+\ .
$$
But $\ufr^+ f = 0$ and $\hfr f = c(\chi)f = mf$. Consequently
$(\zfr f)(gu) = (\zfr f)(g)$. The map $I$ therefore is well
defined; it is $G$-equivariant and continuous by construction.
Finally, it is nonzero since it is easily seen to restrict on
$\Ind^{B}_{P^{-}_{\rm o}}(w\chi)$ to the nonzero map considered in
the proof of Prop. 5.5

\medskip

{\bf Remarks:}

1. (Compare [Jan] II.2 and II.8.23.) Consider the algebraic
character
$$
\chi (\pmatrix{t_1 & 0\cr 0 & t_2})= t_1^{m_1}t_2^{m_2}
$$
with $m_i \in \Zdss$ of $T$. Then $\chi$ is dominant
(for the opposite Borel subgroup $P^-$) if and only if $c(\chi) =
m_2 - m_1 \leq 0$. In this case the $K$-rational representation of
$GL_2$ of highest weight $\chi$ (w.r.t. $P^-$) can be described as
the algebraic induction ${\rm ind}_P^G (\chi)$ which obviously is
a finite dimensional $G$-invariant subspace of $\Ind_P^G (\chi)$;
in particular, $M_{\chi}$ has the finite dimensional
$D(G,K)$-module quotient ${\rm ind}_P^G ((w\chi)^{-1})$.

2. Consider the case $c(\chi) \in -\Ndss_0$. By Thm. 6.1 and Prop.
6.2 the simple $D(G,K)$-module $M_{\epsilon^{1-c(\chi)}
\chi}$ maps isomorphically onto a simple $D(G,K)$-submodule
$M_{\chi}^{(0)}$ in $M_{\chi}$. Morita shows in [Mor] that the
quotient $M_{\chi}/M_{\chi}^{(0)}$ is dual to the $G$-invariant
closed subspace of ``locally polynomial'' functions in
$\Ind_{P}^{G}(\chi)$ (see [Mor] 2-2 for the definition of this
subspace.) Because the space of locally polynomial functions is a
direct limit of finite dimensional spaces, it acquires  the finest
locally convex topology, and so $M_{\chi}/M_{\chi}^{(0)}$ is its
full algebraic dual. Morita analyzes the space of locally
polynomial functions in [Mor] \S7 and shows that it is either
simple as a $(\gfr\lfr_2 (\dQ), G)$ module or has a two step
filtration with simple $(\gfr\lfr_2 (\dQ), G)$-module subquotients
(the bottom one being the algebraic induction discussed in the
first remark). A fortiori, these subquotients are simple as
$D(G,K)$ modules. The subquotients of the corresponding one or
two-step dual filtration on $M_{\chi}/M_{\chi}^{(0)}$, are not in
general simple modules. We also remark that the map $I$ from the
proof of Prop. 6.2 appears in a less conceptual form already in
[Mor] 6-2 and is attributed there to Casselman.

\vfill\eject

{\bf References}

\parindent=23truept

\ref{[Am1]} Amice, Y.: Interpolation $p$-adique. Bull. Soc. math.
France 92, 117-180 (1964)

\ref{[Am2]} Amice, Y.: Duals. Proc. Conf. on $p$-Adic Analysis,
Nijmegen 1978, pp. 1-15

\ref{[B-GAL]} Bourbaki, N.: Groupes et alg\`ebres de Lie, Chap. 1-3. Paris:
Hermann 1971, 1972

\ref{[B-TVS]} Bourbaki, N.: Topological Vector Spaces.
Berlin-Heidelberg-New York: Springer 1987

\ref{[B-VAR]} Bourbaki, N.: Vari\'et\'es diff\'erentielles et
analytiques. Fascicule de r\'e\-sul\-tats. Paris: Hermann 1967

\ref{[GKPS]} De Grande-De Kimpe, N., Kakol, J., Perez-Garcia, C.,
Schikhof, W: $p$-adic locally convex inductive limits. In $p$-adic
functional analysis, Proc. Int. Conf. Nijmegen 1996 (Eds.
Schikhof, Perez-Garcia, Kakol), Lect. Notes Pure Appl. Math., vol.
192, pp.159-222. New York: M. Dekker 1997

\ref{[DG]} Demazure, M., Gabriel, P.: Groupes Alg\'ebriques.
Amsterdam: North-Holland 1970

\ref{[Fe1]} F\'eaux de Lacroix, C. T.: $p$-adische Distributionen.
Diplomarbeit, K\"oln 1992

\ref{[Fe2]} F\'eaux de Lacroix, C. T.: Einige Resultate \"uber die topologischen
Darstellungen $p$-adischer Liegruppen auf unendlich dimensionalen
Vektorr\"aumen \"uber einem $p$-adischen K\"orper. Thesis, K\"oln
1997, Schriftenreihe Math. Inst. Univ. M\"unster, 3. Serie, Heft 23,
pp. 1-111 (1999)

\ref{[Gru]} Gruson, L.: Th\'eorie de Fredholm $p$-adique. Bull. Soc.
math. France 94, 67-95 (1966)

\ref{[Hel]} Helmer, O.: The elementary divisor theorem for certain
rings without chain conditions. Bull. AMS 49, 225-236 (1943)

\ref{[Jan]} Jantzen, J.C.: Representations of Algebraic Groups.
Academic Press 1987

\ref{[Kom]} Komatsu, H.: Projective and injective limits of weakly
compact sequences of locally convex spaces. J. Math. Soc. Japan 19,
366-383 (1967)

\ref{[Laz]} Lazard, M.: Les z\'eros des fonctions analytiques
d'une variable sur un corps valu\'e complet. Publ. math. IHES 14,
47-75 (1962)

\ref{[Mor]} Morita, Y.: Analytic Representations of $SL_2$ over a
$p$-Adic Number Field, III. In Automorphic Forms and Number Theory,
Adv. Studies Pure Math. 7, pp.185-222. Tokyo: Kinokuniya 1985

\ref{[NFA]} Schneider, P.: Nichtarchimedische Funktionalanalysis.
Course at M\"un-ster (1997)

\ref{[Sch]} Schneider, P.: $p$-adic representation theory. The 1999
Britton Lectures at McMaster University. Available at
www.uni-muenster.de/math/u/ schneider

\ref{[ST]} Schneider, P., Teitelbaum, J.: $p$-adic boundary
values. To appear in Ast\'erisque

\ref{[Ti1]} van Tiel, J.: Espaces localement $K$-convexes I-III.
Indagationes Math. 27, 249-258, 259-272, 273-289 (1965)

\ref{[Ti2]} van Tiel, J.: Ensembles pseudo-polaires dans les
espaces localement $K$-convexes. Indagationes Math. 28, 369-373
(1966)

\bigskip

\parindent=0pt

Peter Schneider\hfill\break Mathematisches Institut\hfill\break
Westf\"alische Wilhelms-Universit\"at M\"unster\hfill\break
Einsteinstr. 62\hfill\break D-48149 M\"unster, Germany\hfill\break
pschnei@math.uni-muenster.de\hfill\break
http://www.uni-muenster.de/math/u/schneider\hfill\break

\noindent
Jeremy Teitelbaum\hfill\break Department of Mathematics,
Statistics, and Computer Science (M/C 249)\hfill\break University
of Illinois at Chicago\hfill\break 851 S. Morgan St.\hfill\break
Chicago, IL 60607, USA\hfill\break jeremy@uic.edu\hfill\break
http://raphael.math.uic.edu/$\sim$jeremy\hfill\break

\end